\documentclass[a4paper,12pt]{amsart}
\usepackage{amsmath,amsfonts,amssymb}
\usepackage[latin2]{inputenc} 
\usepackage{accents}

\def\wt#1{\widetilde{#1}}
\def\ol#1{\overline{#1}}

\def\bbR{{\mathbb{R}}}

\newcommand{\comments}[1]{}
\newcommand{\newc}{\newcommand}

\newcommand{\Ker}{\operatorname{Ker}}
\newcommand{\Proj}{\operatorname{Proj}}

\let\ccdot\cdot
\def\cdot{\hbox to 2.5pt{\hss$\ccdot$\hss}}

\newcommand{\be}{\beta}
\newcommand{\ga}{\gamma}
\newcommand{\de}{\delta}

\newcommand{\ka}{\kappa}
\newcommand{\la}{\lambda}
\newcommand{\om}{\omega}
\renewcommand{\phi}{\varphi}
\newcommand{\ph}{\varphi}

\newcommand{\si}{\sigma}

\newcommand{\Ga}{\Gamma}

\newcommand{\Om}{\Omega}

\newc{\aI}{\mbox{\boldmath{$ I$}}}
\newc{\aR}{\mbox{\boldmath{$ R$}}}
\newc{\aDeR}{\mbox{\boldmath{$ U$}}_B{}^P{}_C{}^Q}
\newc{\al}{\mbox{\boldmath$ \Delta$}}             
\newc{\nda}{\mbox{\boldmath$ \nabla$}}
\newc{\ad}{\mbox{\boldmath$ d$}}
\newc{\da}{\mbox{\boldmath$ \delta$}}
\newc{\aK}{\mbox{\boldmath{$ K$}}}
\newc{\aL}{\mbox{\boldmath{$ L$}}}

\newtheorem{theorem}{Theorem}[section]

\newtheorem{proposition}[theorem]{Proposition}
\newtheorem{corollary}[theorem]{Corollary}

\newcommand{\g}{{\frak g}}

\newcommand{\cA}{{\mathcal A}}

\newcommand{\cC}{{\mathcal C}}

\newcommand{\ce}{{\mathcal E}}
\newcommand{\cE}{{\mathcal E}}

\newcommand{\cQ}{{\mathcal Q}}

\newcommand{\cT}{{\mathcal T}}

\newcommand{\Ps}{\Psi}
\newcommand{\Rho}{{\mbox{\sf P}}}
\newcommand{\Up}{\Upsilon}

\newcommand{\End}{\operatorname{End}}
\newcommand{\im}{\operatorname{im}}
\newcommand{\Ric}{\operatorname{Ric}}

\newcommand{\wh}{\widehat}
\newcommand{\R}{\mathbb{R}}

\newcommand{\vol}{\mbox{\large\boldmath $ \epsilon$}}

\newcommand{\nn}[1]{(\ref{#1})}

\newcommand{\X}{\mbox{\boldmath{$ X$}}}

\newcommand{\Y}{\mbox{\boldmath{$ Y$}}}

\newcommand{\bg}{\mbox{\boldmath{$ g$}}}

\newcommand{\V}{{\mbox{\sf P}}}                   
\newcommand{\J}{{\mbox{\sf J}}}

\newc{\strutdd}{\rule{0mm}{5mm}}

\newcommand{\lpl}                         
{\mbox{$
\begin{picture}(12.7,8)(-.5,-1)
\put(2,0.2){$+$}
\put(6.2,2.8){\oval(8,8)[l]}
\end{picture}$}}

\usepackage{ifthen}

\newc{\tensor}[1]{#1}

\newc{\Mvariable}[1]{\mbox{#1}}

\newc{\down}[1]{{}_{
\ifthenelse{\equal{#1}{;}}{|}{#1}}}

\newc{\up}[1]{{}^{#1}}
\newc{\C}{C}

\newc{\JulyStrut}{\rule{0mm}{6mm}}
\newc{\midtenPan}{\mbox{\sf S}}
\newc{\midten}{\mbox{\sf T}}
\newc{\midtenEi}{\mbox{\sf U}}
\newc{\ATen}{\mbox{\sf E}}
\newc{\BTen}{\mbox{\sf F}}
\newc{\CTen}{\mbox{\sf G}}

\def\sideremark#1{\ifvmode\leavevmode\fi\vadjust{\vbox to0pt{\vss
 \hbox to 0pt{\hskip\hsize\hskip1em
 \vbox{\hsize3cm\tiny\raggedright\pretolerance10000
 \noindent #1\hfill}\hss}\vbox to8pt{\vfil}\vss}}}%

\input xy
\xyoption{all}

\newtheorem{definition}[theorem]{Definition}

\theoremstyle{definition}
\newtheorem{example}[theorem]{Example}
\newtheorem{remark}[theorem]{Remark}

\newcommand{\mV}{\mathbb V}

\newcommand{\bd}{\begin{definition}}
\newcommand{\ed}{\end{definition}}

\newcommand{\br}{\begin{remark}}
\newcommand{\er}{\end{remark}}
\newcommand{\gog}{{\mathfrak g}}

\newcommand{\bt}{\begin{tabular}}
\newcommand{\et}{\end{tabular}}

\renewcommand{\ker}{\operatorname{ker}}

\renewcommand{\Im}{\operatorname{Im}}

\newcommand{\Id}{\operatorname{id}}

\def\cdot{\hbox to 2.5pt{\hss$\ccdot$\hss}}

\def\gog{\frak{g}}

\newcommand{\pa}{\partial}

\newcommand{\io}{\iota}
 
\def\C{\mathbb{C}}

\def\R{\mathbb{R}}
\def\V{\mathbb{V}}
\def\W{\mathbb{W}}
\def\X{\mathbb{X}}
\def\Y{\mathbb{Y}}
\def\Z{\mathbb{Z}}

\def\cC{\mathcal{C}}

\def\cE{\mathcal{E}}
\def\cG{\mathcal{G}}

\def\cQ{\mathcal{Q}}

\def\cT{\mathcal{T}}

\def\al{\alpha}
\def\be{\beta}
\def\ga{\gamma}
\def\de{\delta}

\def\et{\eta}

\def\io{\iota}
\def\ka{\kappa}
\def\la{\lambda}
\def\rh{\rho}
\def\si{\sigma}

\def\ph{\varphi}

\def\om{\omega}
\def\Ga{\Gamma}

\def\Ps{\Psi}
\def\Om{\Omega}
\def\Up{\Upsilon}
\def\na{\nabla}
\def\pa{\partial}

\def\form#1{\mathbf{#1}}
\def\dform#1{\dot{\mathbf{#1}}}
\def\ddform#1{\ddot{\mathbf{#1}}}
\def\dddform#1{\dddot{\mathbf{#1}}}

\def\sideremark#1{\ifvmode\leavevmode\fi\vadjust{\vbox to0pt{\vss
 \hbox to 0pt{\hskip\hsize\hskip1em
 \vbox{\hsize3cm\tiny\raggedright\pretolerance10000
  \noindent #1\hfill}\hss}\vbox to8pt{\vfil}\vss}}}%

\def\idx#1{{\em #1\/}}

\author{M.\ Hammerl, P.\ Somberg, V.\ Sou\v cek, J.\ \v Silhan}
\title{Invariant prolongation of overdetermined PDE's in projective, 
conformal and Grassmannian geometry}
\date{\today}

\begin{document}

\maketitle

\begin{abstract}
  This is the second in a series of papers on natural modification of the
  normal tractor connection in a parabolic geometry, which naturally prolongs
  an underlying overdetermined system of invariant differential equations. 
  We give
  a short review of the general procedure developed in \cite{hsss} and then 
  compute the prolongation covariant derivatives 
  for a number of interesting examples in projective, conformal and
  Grassmannian geometries.
\end{abstract}

\section{Introduction}

In this paper we study certain overdetermined linear systems of PDE's
that have geometric origin and satisfy strong invariance properties.
The goal is to rewrite these systems in a closed form, which for our purposes
means to find an extended system described by a covariant derivative in such a way 
that parallel sections with respect to this covariant derivative are
in one to one correspondence with solutions of the original equation.
The main advantage of such a prolongation is clear - one immediately
obtains a bound on the dimension of the solution space and the curvature
of this covariant derivative obstructs the existence
of a solution. Moreover, there is a neat relationship between geometry
of the underlying manifold and the extended prolongation system, see e.g.
\cite{prolong},\cite{hsss} and the references therein.

The equations we study appear naturally for parabolic geometries like
projective, conformal or Grassmannian structures and include as a special instances
the equations describing the infinitesimal symmetries of
geometric structures. Special examples of overdetermined linear systems
of invariant equations coming from parabolic geometries are discussed in e.g.,
\cite{prolong}, \cite{EaMa}, \cite{inaut},  \cite{DuTo4dim},  \cite{GoSi}.

In fact, the invariant equations in question appear in the Bernstein-Gelfand-Gelfand 
(BGG for short) sequences, which are the source of overdetermined invariant operators 
resp.\ their prolonged systems in question.
The prolongation of the first operator in the BGG sequence is realized 
by certain commutative square related to BGG operators in the sequence.
We are constructing also examples of commutative squares for all operators 
in the BGG sequence.

\subsection{The BGG-sequence}
Let $G$\ be a semi-simple Lie group and $P\subset G$\ a parabolic subgroup.
A \emph{parabolic geometry} on a manifold $M$\ consists of a
$P$-principal bundle $\mathcal{G}\to M$\ together with a \emph{Cartan
connection $1$-form} $\om\in\Om^1(\mathcal{G},\g)$, \cite{CSbook}.
Here $\frak{g}$ denotes the Lie algebra of $G$.
A major development in the construction of differential invariants of 
parabolic structure was done in \cite{BGG-2001}, and the construction 
was subsequently simplified in \cite{BGG-Calderbank-Diemer}.

\def\dcov{{d^{\nabla}}}
\def\delstar{\pa^*}
Let $\mathbb{V}$\ be a finite dimensional $G$-representation. It is well known that
the associated \emph{tractor bundle}  $V=\mathcal{G}\times_P\mathbb{V}$\ carries
the canonical \emph{tractor covariant derivative} $\na$ 
induced by the Cartan connection form $\om$, see e.g. \cite{thomass}.
The connection uniquely extends to an exterior covariant derivative
on the spaces $\ce^k(V):=\Om^k(M,V)$\ of $k$-forms with values in the vector bundle $V$, 
denoted $\dcov:\ce^k(V)\to\ce^{k+1}(V)$.
The lowest homogenoues part of $\dcov$ is the $G_0$-equivariant Lie algebraic
differential $\pa_k: \ce^k(V)\to\ce^{k+1}(V)$ termed the 
\emph{Kostant differential}, \cite{Kostant}.
Here $G_0$ denotes the Levi part of $P$.
Its adjoint, the \emph{Kostant codifferential} $\pa_k^*$ is $P$-equivariant 
and gives rise to a complex
\begin{align*}
  \ce^{k+1}(V)\overset{\delstar_{k+1}}{\to}\ce^k(V),\,\, \delstar_k\circ\delstar_{k+1}=0.
\end{align*}
There are Lie algebra cohomology bundles $H_k=\ker\delstar_k/\im\delstar_{k+1}$ due to 
the $P$-equivariant projection
\begin{align*}
  \Pi_k:\ker\delstar_k\to H_k.
\end{align*}
The basic ingredient of the BGG-machinery are the differential \emph{BGG-splitting operators}
\begin{align*}
  L_k:H_k\to\ker\delstar_k,
\end{align*}
defined uniquely by the property that for every smooth section $\si\in \Ga(H_k)$\ one
has 
$$
\delstar_{k+1}(d^\nabla (L_k(\si)))=0.
$$ 
In particular, one can form the \emph{BGG-operators}
\begin{align*}
  &D_k:H_k\to H_{k+1},\ D_k:=\Pi_{k+1}\circ \dcov\circ L_k.
\end{align*}
It will be usually clear from the context what is the appropriate
value for homogeneity $k$ of the form which is acted upon by any 
of operators, i.e. we usually omit this subscript from the notation. 

Let us briefly review the invariant prolongation procedure obtained in \cite{hsss}:

\subsection{Prolongation of the first BGG operator $D_0$}\label{prolongD0}
The \emph{first BGG-operator}\ $D_0$\ associated to $\mathbb{V}$ is overdetermined,
and our aim is the construction of invariant prolongation of the corresponding systems
$D_0\si=0$\ on $\si\in \Gamma (H_0)$. 
Let us recall that the approach of \cite{hsss} starts by introducing certain
class of linear connections on $V$ which are modifications of tractor covariant
derivative
$\na^V$. The first condition on a modification map $\Phi\in\cE^1(\End V)$\
is that it is homogeneous of degree $\geq 1$ with respect to the natural
filtrations on $TM$\ and $V$, for which we write $\Phi\in(\cE^1(\End V))^1$.
This ensures that basic constructions of the BGG-machinery still work.
The next condition is that for any section $s\in\Ga(V)$\ we have
that $\Phi s\in\ce^1(V)$\ has values in $\im\delstar$. As a consequence,
the modified covariant derivative is in a suitable sense compatible
with the underlying first BGG-operator $D_0$. The latter condition
can be rewritten as $\Phi\in\Im(\pa^*_V\otimes \Id_{V^*})$, 
thus we arrive at a class of \emph{admissible} covariant derivatives
\begin{align*}
  \cC=
\left\{
\wt{\na} = \na+\Phi|\Phi\in Im \; (\pa^*_V\otimes \Id_{V^*})\cap (\cE^1(\End V))^1%
\right\}.
\end{align*}
Here $\pa^*_V$ denotes $\pa^*$ acting on $\cE^1(V)$ (and \idx{not} on 
$\cE^1(\End V)$) and the same applies for $\pa^*_V$ acting on $\cE^k(V)$.

\def\ti{\tilde}
The main theorem of \cite{hsss} is then
\begin{theorem}\label{thm-prolongcon}
There exists a unique covariant derivative $\ti\na\in\cC$ characterized by the property
$$
(\pa^*_V\otimes \Id_{V^*})(\wt{\Om})=0,
$$
where $\wt{\Om}$ is the curvature of $\ti\na$. 
\end{theorem}
This implies
$\ti\na\circ L_0=L_1\circ D_0$,
which in turn yields
\begin{corollary}\label{cor-prolong}
Consider a tractor bundle $V$ and the covariant derivative $\ti\na$\ in Theorem
\ref{thm-prolongcon}.
Then $\ti\na$ gives a prolongation of the first BGG operator $D_0$ in the 
sense that
the restriction of the projection $\Pi_0:V\to H_0$ to $\ti\na$-parallel
sections is an isomorphism with the kernel of $D_0$ acting on smooth sections $\Gamma (H_0)$ and inverted
by the differential splitting operator $L_0:H_0\to V$.
\end{corollary}
We therefore say that $\ti\na$\ is the \emph{prolongation covariant derivative}.

\subsection{Commutativity for all $D_k$}\label{prolongDk}
In \cite{hsss} the authors also obtained the analogue of $\wt{\na}$ on 
$\cE^k(V)$. Here $d^\na$ gives rise to the class 
$$
\cC_k:=\{ \tilde{d}_k=d^\na+\Phi \mid \Phi\in A^1, \Im \Phi\subset \Im \pa^*\}
$$
where $A:=Hom(\cE^k(V),\cE^{k+1}(V))$ and $A^1$ denotes homomorphisms
homogeneous of the degree $\geq 1$. Then it turns out
there is a unique $\tilde{d}_k \in \cC_k$ such that 
$\pa_V^* \circ d^\na \circ \tilde{d}_k =0$. This then implies
\begin{align*}
  \ti d_k \circ L_k=L_{k+1}\circ D_k
\end{align*}
and $\Pi_k$\ and $L_k$\ restrict to inverse isomorphisms between 
$\Ker \ti d_k\cap\Ker\pa^*$  and  $\Ker D_k$.

\subsection{The guideline for computing examples}
Here is the manual for treating particular examples, which can 
be used to derive the explicit form of the prolongation covariant derivative.
In practice, the normalization procedure for canonical tractor covariant derivative 
can be summarized as an algorithm based on the following list of steps:
\begin{enumerate} 
\item
Choose a parabolic geometry $(\cG, P, M, \omega)$, where $\cG\to M$ is 
a principal $P$-bundle on $M$ and $\omega\in\Omega^1(\cG ,\gog)$.
Choose also a finite dimensional $G$-module $\mV$ and its associated vector
bundle $V$ termed tractor bundle. Let us fix the two consecutive vector 
bundles of $k$ resp. $(k+1)$-forms twisted by $V$. 
\item
Decompose both spaces of $k$ resp. $(k+1)$-forms twisted by $V$ with respect to
$G_0$, the Levi factor of the parabolic subgroup $P$. Then compute the value 
of the Laplace-Kostant algebraic operator $\square$ associated to $\delstar$ 
on each irreducible $G_0$-summand 
(i.e.\ $G_0$-graded components associated to $P$-equivariant filtration) 
either by evaluating the action of Casimir operator or 
from the definition $\square = \pa^* \pa + \pa \pa^*$.

\item  Choose a Weyl structure, so that there is a well defined
 splitting of the filtered bundle $V$ into a direct sum of homogeneous components.

\item Now the procedure splits into two cases:

\begin{itemize}
\item
The computation of the prolongation covariant derivative.
 
Check, if  $(\pa_V^* \otimes \Id_{V^*})(\Om)$, where $\Om$ is the curvature of $\na$, is trivial. 
In positive case, the procedure
ends and we have computed the prolongation covariant derivative. 

If $\al:=(\pa_V^*\otimes \Id_{V^*})(\Om)\not=0,$  take the lowest nontrivial homogeneous part
$\al_j$ of $\al$ and define
$$
\Phi=-\square^{-1}\al_j;\;\; \na'=\na+\Phi.
$$

Then repeat the procedure with $\na$ replaced by $\na'.$ By construction, the lowest nontrivial
component of $\al$ in the next step will have degree higher then in the previous step, hence the procedure will terminate in a finite number of steps (bounded by the length of the grading
of $\V$).

\item The case of the whole sequence of commuting squares.

Here we use another procedure based on the following algorithm.
Consider two consecutive squares containing the exterior covariant derivatives
 $d^\na_k:\cE^k(V)\mapsto \cE^{k+1}(V)$ and $d^\na_{k+1}:\cE^{k+1}(V)\mapsto \cE^{k+2}(V).$ 
 First check, if 
 $$
 (\pa^*_V\otimes \Id_{V^*})(d^\na_{k+1}\circ d^\na_k)
 $$ 
 is trivial.
 If not, the first step is the same as for the construction of prolongation covariant
 derivative above. 
 Consider $\al:=(\pa^*\otimes \Id_{V^*})(d^\na_{k+1}\circ d^\na_k)\not=0,$  take the lowest nontrivial homogeneous part
$\al_j$ of $\al$ and define
$$
\Phi=-\square^{-1}\al_j;\;\; d_k'=d_k^\na+\Phi.
$$

If $\al':=(\pa^*\otimes \Id_{V^*})(d_k^\na\circ d_k')$ is trivial, the procedure terminates and
we define $\tilde{d}_k=d_k'.$ If not, take the lowest nontrivial homogeneous part
$\al'_{j'}$ of $\al'$ and define
$$
\Phi'=-\square^{-1}\al'_{j'};\;\; d_k''=d_k'+\Phi'.
$$
By construction, the degree $j'$ will be bigger than $j,$ 
hence the procedure will terminate in a finite number of steps (bounded again by the length of the grading of $V$).
Note that iterations $\phi$ here are, in general, differential operators and their
order rises (in general) by one with each iteration.
\end{itemize}
\end{enumerate}

The panorama of examples presented in this article follows criterions to be
useful, nonelementary, going beyond the examples scattered in the references  
and at the same time computable by hand while demonstrating the powerful 
machine developed in \cite{hsss}.
The interested reader will easily recognize the complexity of the 
computation both in general and specific situations of interest.

\section{Notation}

In this section we review the basic notation and conventions
related to the results of our article. 

\subsection{Forms, tensors and tensorial actions} \label{forms}
In order to be
explicit and efficient in calculations involving bundles of possibly
high rank it is necessary to introduce some further abstract index
notation.  In the usual abstract index conventions one would write
$\cE_{[ab\cdots c]}$ (where there are implicitly $k$-indices skewed
over) for the space $\cE^k$. To simplify subsequent expressions we
 use the following conventions. Firstly
indices labeled with sequential superscripts which are
at the same level (i.e. all contravariant or all
covariant) indicate a completely skew set of indices.
Formally we set $a^1 \cdots a^k = [a^1 \cdots a^k]$ and so, for example,
$\cE_{a^1 \cdots a^k}$ is an alternative notation for $\cE^k$
while $\cE_{a^1 \cdots a^{k-1}}$ and $\cE_{a^2 \cdots a^k}$ both denote 
$\cE^{k-1}$. Next we
abbreviate this notation via multi-indices: We will use the form
indices
$$
\begin{aligned}
\form{a}^k &:=a^1 \cdots a^k =[a^1 \cdots a^k], \quad k \geq 0,\\
\dot{\form{a}}^k &:= a^2 \cdots a^k=[a^2 \cdots a^k], \quad k \geq 1,\\
\ddot{\form{a}}^k &:= a^3 \cdots a^k =[a^3 \cdots a^k], \quad k \geq 2,  \\
\dddot{\form{a}}^k &:= a^4 \cdots a^k =[a^4 \cdots a^k], \quad k \geq 3.  \\
\end{aligned}
$$ If, for example, $k=1$ then $ \dot{\form{a}}^k$ simply means the
index is absent, whereas if $k=1$ then $\ddform{a}$ means the term
containing the index $\ddform{a}$ is absent. For
example, a 3--form $\ph$ can have the following possible equivalent
structures of indices:
$$ \ph_{a^1a^2a^3} = \ph_{[a^1a^2a^3]} = \ph_{\form{a}^3}
= \ph_{a^1\dot{\form{a}}^3}
= \ph_{[a^1\dot{\form{a}}^3]} =
   \ph_{a^1a^2\ddot{\form{a}}^3} \in \mathcal{E}_{\form{a}^3} = \mathcal{E}^3.
$$
Note the exterior derivative $d$ on a $k$-form $f_\form{a}$ can be written 
as $(df)_{a^0\form{a}} = \na_{a^0} f_\form{a}$
for any torsion--free affine connection $\na$.

Later on we define the standard tractor bundle denoted by $\cE^A$ and 
its dual $\cE_B$. The form index notation developed above will be used also
for skew symmetric powers of these bundles. For example, \ the bundle of
\idx{tractor $k$--forms} $\cE_{[A^1\cdots A^k]}$ will be
denoted by $\cE_{A^1\cdots A^k}$ or $\mathcal{E}_{\form{A}^k}$.

The bundle of endomorphisms of $\cE^A$ (or $\cE_A$), $\cE^E{}_F$, 
clearly injects $\cE^E{}_F \subseteq \End(\cT)$ for any tractor bundle 
$\cT \subseteq (\bigotimes \cE^A) \otimes ( \bigotimes \cE_B)$.
Consider $\ga^E{}_F \in \cE^E{}_F$ and $f \in \cT$. The endomorphism
$\ga$ acts on $\cT$ and we denote this action by $\sharp$. 
That is, $\ga \sharp f \in \cT$. Using the abstract tractor indices,
$\sharp$ is given by the usual tensorial action, i.e.\
$(\ga \sharp f)^A = \ga^A{}_P f^P$ for $f^A \in \cE^A$ and
$(\ga \sharp f)_A = -\ga^P{}_A f_P$ for $f_A \in \cE_A$.
One then computes $\sharp$ on the tensor products of $\cE^A$ and 
$\cE_B$ using the Leibniz rule.
We further put $\ga \sharp$ to be zero on $\cE^a$, $\cE_b$ and density 
bundles (which we introduce later) and, using the Leibniz rule,
extend $\ga \sharp$ to the tensor
products of $\cT$ with latter three bundles. 
Finally note the action $\sharp$ is denoted $\bullet$ in 
\cite{CSbook}.

\subsection{The adjoint tractor bundle and the Laplace-Kostant operator}
\label{KoLap}

The bundle $\cA = \cG \times_P \frak{g}$ is called the \idx{adjoint} tractor 
bundle. By definition,  
$\cA \subseteq \cE^A{}_B$ and more generally  
$\cA \hookrightarrow \End(\cT)$ for any tractor bundle $\cT$.
We shall use $\sharp$ to denote the action of sections of $\cA$ on 
$\cT$ as introduced above. 
Note the curvature of the normal tractor covariant derivative $\na$
is the section of $\cE_{a^0a^1} \otimes \cA$ and the curvature action 
is $2 (d^\na \na f)_{a^0a^1} 
= 2\na_{a^0} \na_{a^1} f  = (\Om \sharp f)_{a^0a^1} \in \cE_{[ab]} \otimes \cT$ 
for each $f \in \cT$.

We have identifications
$\cE_a \cong \cG \times_P \frak{g}_-$ and 
$\cE^a \cong \cA/\cA'$, $\cA' := \cG \times_P \frak{p}$,
which allow to define inclusions
$\io: \cE_a \hookrightarrow \cA$ and 
$\bar{\io}: \cE^a \hookrightarrow \cA/\cA'$. (The latter is just the identity.)
We extend these inclusions to 
$$
\io: \cE_\form{a} \hookrightarrow \cE_{\dform{a}} \otimes \cA
\quad \text{and} \quad
\bar{\io}: \cE_\form{a} \stackrel{\textstyle \de_{a^0}{}^b}{\longrightarrow} 
\cE_{a^0\form{a}}{}^b \hookrightarrow \cE_{a^0\form{a}} \otimes \cA/\cA'.
$$
Recall that here and below, we use a chosen Weyl structure and the 
corresponding splittings.

Our aim is to use these tools to express
Kostant's differential $\pa$, codifferential $\pa^*$ 
and in particular the Laplace-Kostant operator $\Box$ \cite{Kostant} 
in a form suitable for computations in abstract indices.
Defined on $\cE_\form{a} \otimes \cT$,
$\form{a} = \form{a}^k$ for any tractor bundle $\cT$, they have the form 
\begin{align*}
&\pa: \cE_\form{a} \otimes \cT \stackrel{\bar{\io}}{\hookrightarrow}
\cE_{a^0\form{a}} \otimes \cA/\cA' \otimes \cT 
\stackrel{\sharp}{\longrightarrow} 
\cE_{a^0\form{a}} \otimes \cT, \\
&\pa^*: \cE_\form{a} \otimes \cT \stackrel{\io}{\hookrightarrow}
\cE_{\dform{a}} \otimes \cA \otimes \cT \stackrel{\sharp}{\longrightarrow} 
\cE_{\dform{a}} \otimes \cT \quad \text{and} \\
&\Box_k = \pa \pa^* + \pa^* \pa: \cE_\form{a} \otimes \cT \longrightarrow
\cE_{\form{a}} \otimes \cT.
\end{align*}
Note $\pa^*$ is invariant but $\pa$ (thus also $\Box_k$) depends on the 
choice of splitting of the tractor bundles in question. 
However, $\Box_k$ is invariant on completely reducible subquotients of 
$\cE_\form{a} \otimes \cT$ and acts by a scalar multiple on each 
irreducible component of such subquotients. That is, we choose a splitting
of the tractor bundle $\cE_\form{a} \otimes \cT$ to compute $\Box_k$ but
the value of $\Box_k$ on a given completely reducible subquotient alone
is independent of this choice.

The symbol $\lpl$ denotes the composition $P$-module structure of 
representations or vector bundles.

Finally note one can compute $\Box_k$ from highest weight of bundles 
concerned, see \cite{Kostant}. We shall use this (less explicit) approach in cases
when the abstract index computation is getting too complicated.

\vspace{1ex}

Now we are ready to discuss specific geometries. In each case, we first
summarize the tractor calculus. We shall particularly need the normal 
tractor covariant derivative $\na$ and the Kostant's differential and 
codifferential $\pa$ and $\pa^*$, respectively. Using these we compute
the prolongation covariant derivative
$\wt{\na}$ and/or $\tilde{d}$ on certain bundles.

\section{Projective geometry}

We follow the notation from \cite{thomass} here. The projective 
structure on a smooth manifold $M$ is given by a class $[\na]$ of 
projectively equivalent torsion free connections. That is, 
connections $\hat{\na} \in [\na]$ are parametrised by one forms
$\Up_a \in \cE_a \cong \Ga(T^*M)$ and have the form
\begin{align}
&\hat{\na}_a \ph = \na_a \ph + w\Up_a \ph, 
\quad \ph \in \cE(w), \notag \\
&\hat{\na}_a f^b = \na_a f^b + \Up_a f^b + \Up_c f^c \de_a^b,
\quad f^b \in \cE^b \label{trans} \\
&\hat{\na}_a \om_b = \na_a \om_b - \Up_a \om_b - \Up_b \om_a, 
\quad \om_a \in \cE_a. \notag
\end{align}

The curvature tensor $R_{ab}{}^c{}_d$ of a torsion free $\na$ is defined by 
$(\na_a \na_b - \na_b \na_a) f^c = R_{ab}{}^c{}_p f^p$ and it decomposes
$$
R_{ab}{}^c{}_d = W_{ab}{}^c{}_d + 2\de_{[a}{}^c \Rho_{b]d} + \be_{ab}\de^c{}_d,
\quad \be_{ab} = -2\Rho_{[ab]}.
$$
Here $W_{ab}{}^c{}_d$ is projectively invariant (and irreducible) Weyl 
tensor,  $\Rho$ is the Schouten tensor, 
$\hat{\Rho}_{ab} = \Rho_{ab} - \na_a \Up_b  + \Up_a\Up_b$ and 
$\hat{\be}_{ab} = \be_{ab} + 2\na_{[a}\Up_{b]}$.
We put $A_{abc} := 2\na_{[a}\Rho_{b]c}$. Then
the Bianchi identity $\na_{[a}R_{bc]}{}^d{}_e=0$ implies
$$
\na_{c}W_{ab}{}^c{}_d = (n-2) A_{abd}
\quad \text{and} \quad \na_{[a} \be_{cd]}=0.
$$
The cohomology class $[\be] \in H^2(M,\bbR)$ is a global invariant of 
the projective structure. Moreover, 
$(\na_a \na_b - \na_b \na_a) \ph = w \be_{ab} \ph$ for $\ph \in \cE(w)$.

\subsection{Projective tractors}
We shall write sections of the standard projective tractor bundle 
$\cE^A = \cE^a[-1] \lpl \cE[-1]$, resp.\ its dual 
$\cE_A = \cE[1] \lpl \cE_a[1]$ using the injectors $Y^A$, $X^A$,
resp.\ $Y_A$, $X_A$ as 
$$
\begin{pmatrix} 
\si^a \\ \rh
\end{pmatrix}
= Y^A_a \si^a + X^A \rh \in \cE^A,
\quad \text{resp.} \quad
\begin{pmatrix} 
\nu \\ \mu_a
\end{pmatrix}
= Y_A \nu + X_A^a \mu_a \in \cE_A.
$$
Such splittings of $\cE^A$ and $\cE_A$ are parametrised by choices of projective
connections and we call them \idx{projective splittings}.
The change of the splitting under change of the connection parametrised by
$\Up_a \in \cE_a$ is
\begin{align*}
&\widehat{
\begin{pmatrix} 
\si^a \\ \rh
\end{pmatrix}
} = 
\begin{pmatrix} 
\si^a \\ \rh - \Up_a \si^a  
\end{pmatrix},
\ \ \mbox{i.e.} \ \
\hat{Y}^A_a = Y^A_a + X^A \Up_a,\ 
\hat{X}^A = X^A
\quad \text{and} \\
&\widehat{
\begin{pmatrix} 
\nu \\ \mu_a
\end{pmatrix}
} = 
\begin{pmatrix} 
\nu \\ \mu_a + \Up_a \nu  
\end{pmatrix},
\ \ \mbox{i.e.} \ \
\hat{Y}_A = Y_A - X_A^a \Up_a,\ 
\hat{X}_A^a = X_A^a.
\end{align*}
That is,  $X^A \in \cE^A[1]$, $X_A^a \in \cE_A^a[-1]$ are invariant
and $Y^A_a \in \cE^A_a[1]$, $Y_A \in \cE_A[-1]$ depend on the choice
of the projective scale. 
We assume the normalisation of these such that 
$Y_A X^B + X_A^c Y^B_c = \de_A{}^B$, i.e.\ $Y_C X^C = 1$ and 
$X_C^a Y^C_b = \de^a{}_b$.

The normal covariant derivative is given by
\begin{align*}
&\na_c
\begin{pmatrix} 
\si^a \\ \rh 
\end{pmatrix} = 
\begin{pmatrix} 
\na_c \si^a + \rh \de_c{}^a \\ \na_c \rh - P_{cp} \si^p
\end{pmatrix} 
\quad \text{and} \quad
\na_c 
\begin{pmatrix} 
\nu \\ \mu_a 
\end{pmatrix} = 
\begin{pmatrix} 
\na_c \nu - \mu_c  \\ \na_c \mu_a + P_{ca} \nu
\end{pmatrix}, \quad \text{i.e.}\\
&\na_c^{} Y^A_a = -X^A \Rho_{ca},\ \na_c X^A = Y^A_c
\ \ \text{and} \ \
\na_c^{} Y_A = X_A^a \Rho_{ca},\ \na_c X_A^a = -Y^A \de_c^a.
\end{align*}
and its $\Om$ curvature has the form
$$
\Om_{ab}{}^E{}_F = Y^E_e X_F^f W_{ab}{}^e{}_f - X^E X_F^f A_{abf} \in
\cE_{[ab]} \otimes \cA.
$$
That is, $\cA = \text{trace-free} (\cE^E{}_F)$ is the projective adjoint 
tractor bundle where ``trace-free'' denotes the trace--free part.
Hence the curvature action on $\cE_C$ is
$(\na_a\na_b - \na_b\na_a) F_C = (\Om \sharp F)_{abC} = 
- \Om_{ab}{}^D{}_C F_D$. We shall often write
$\Om_{ab} \sharp F_C$ instead  of $(\Om \sharp F)_{abC}$ to simplify
the notation. 

\vspace{1ex}

Using the notation developed above, the inclusions $\io$ and $\bar{\io}$
defined in \ref{KoLap} have the form
$Y^E_{a^0} Y_F^{}: \cE_\form{a} \stackrel{\bar{\io}}{\to} 
\cE_{a^0\form{a}}{}^E{}_F$ and 
$X^E X_F^{a^1}: \cE_\form{a} \stackrel{\io}{\to} \cE_{\dform{a}}{}^E{}_F$.
Thus
\begin{align*}
&\pa: \cE_\form{a} \otimes \cT \ni f_\form{a} \mapsto
Y^E_{a^0} Y_F^{} f_\form{a} \stackrel{\sharp}{\longrightarrow} 
\cE_{a^0\dform{a}} \otimes \cT \quad \text{and} \\
&\pa^*: \cE_\form{a} \otimes \cT \ni f_\form{a} \mapsto
X^E X_F^{a^1} f_\form{a} \stackrel{\sharp}{\longrightarrow} 
\cE_{\dform{a}} \otimes \cT
\end{align*}
and we can easily compute $\Box_k$ on $\cE_\form{a} \otimes \cT$ 
using the action $\sharp$
as demonstrated by the following example.

\begin{example} \label{standard}
We shall compute the case $\cT = \cE^C$ in details. Then
$\cE_\form{a}{}^C = \cE_\form{a}{}^c[-1] \lpl \cE_\form{a}[-1]$,
where $\cE_\form{a}$ is irreducible and $\cE_\form{a}{}^c$ has two
irrreducible components (the trace and trace--free parts). 
We shall compute $\Box_k$ separately for all three irreducible components.

We start with (not necessarily irreducible) section 
$\si_\form{a}{}^c \in \cE_\form{a}{}^c[-1]$. Then $\pa$ on
$f_\form{a}{}^C := Y^C_{\,c} \si_\form{a}^{}{}^c$ is zero and 
$X^E X_F^{a^1} \sharp Y^C_{\,c} \si_\form{a}^{}{}^c = 
X^C \si_{p\dform{a}}{}^p = (\pa^* f)_{\dform{a}}{}^C$.
Thus $\pa^*f=0$ for trace--free section $\si_\form{a}{}^c$. 
Assume $\si_\form{a}{}^c = \de_{a^1}^c \tilde{\si}_{\dform{a}}^{}$. Then
$f_\form{a}{}^C = Y^C_{a^1} \tilde{\si}_{\dform{a}}^{}$, 
$(\pa^*f)_{\dform{a}}{}^C = \frac{n-k+1}{k} X^C \tilde{\si}_{\dform{a}}$
thus 
$(\Box_kf)_\form{a}{}^C = (\pa\pa^* f)_\form{a}{}^C = 
Y^C_{a^1} \tilde{\si}_{\dform{a}}^{}$.
Finally if $\bar{f}_\form{a}{}^C = X^C \rh_\form{a}$ then 
$(\pa^*\bar{f})_{\dform{a}}{}^C=0$,
$(\pa \bar{f})_\form{a}{}^C = Y^C_{a^0} \rh_\form{a}^{}$ and
$(\Box_k \bar{f})_\form{a}{}^C = (\pa^*\pa \bar{f})_{\form{a}}{}^C = 
\frac{n-k}{k+1} X^C \rh_\form{a}$

Summarizing, $\Box_k$ acts by zero on the trace--free part of 
$\cE_\form{a}{}^c[-1] = \cE_\form{a}{}^C/\cE_\form{a}[-1]$, by 
$\frac{n-k+1}{k}$ on the trace part, i.e.\
on $\cE_{\dform{a}}[-1] \subseteq \cE_\form{a}{}^C/\cE_\form{a}[-1]$ and by 
$\frac{n-k}{k+1}$ on $\cE_\form{a}[-1] \subseteq \cE_\form{a}{}^C$. 
Note the inclusion $\cE_\form{a}[-1] \hookrightarrow \cE_\form{a}{}^C$
is realized by $X^C: \cE_\form{a}[-1] \to \cE_\form{a}{}^C$.
\end{example}

\subsection{Skew symmetric tractors and tractor forms}

The notation for the standard tractor bundle $\cE^C$ developed above can be 
easily generalised to the products $\bigwedge^\ell \cE^C = \cE^\form{C} =
\cE^\form{c}(-\ell) \lpl \cE^{\dform{c}}(-\ell)$, where
$\form{C} = \form{C}^\ell$. Note 
$\bigwedge^\ell \cE^C \cong \bigwedge^{n-\ell+1} \cE_D$, hence these products
are isomorphic to tractor forms. We put
$$
\Y^\form{C}_{\,\form{c}} = Y^{[C^1}_{\,\ c^1} \ldots Y^{C^\ell]}_{\,c^\ell}
\in \cE^\form{C}_{\,\form{c}}(\ell), \quad
\X^\form{C}_{\,\dform{c}} = X^{[C^1}_{} Y^{C^2}_{\,c^2} \ldots 
Y^{C^\ell]}_{\,c^\ell} 
\in \cE^\form{C}_{\,\dform{c}}(\ell),
$$
and write the sections of $\cE^\form{C}$ as
$$
\begin{pmatrix} 
\si^\form{c} \\ \rh^{\dform{c}}
\end{pmatrix}
= \Y^\form{C}_{\,\form{c}} \si^\form{c} 
+ \X^\form{C}_{\,\dform{c}} \rh^{\dform{c}} \in \cE^\form{C}, \quad
\si^\form{c} \in \cE^\form{c}(-\ell),\ \ 
\rh^{\dform{c}} \in \cE^{\dform{c}}(-\ell)
$$
where $\form{c} = \form{c}^\ell$. The change of the projective rescaling 
parametrised by $\Up_a$ is 
\begin{align*}
&\widehat{
\begin{pmatrix} 
\si^\form{c} \\ \rh^{\dform{c}}
\end{pmatrix}
} = 
\begin{pmatrix} 
\si^\form{c} \\ \rh^{\dform{c}} - \ell \Up_{c^1} \si^\form{c}  
\end{pmatrix},
\ \ \mbox{i.e.} \ \
\hat{\Y}^\form{C}_{\,\form{c}} = \Y^\form{C}_{\,\form{c}} 
+ \ell\, \Up_{c^1} \X^\form{C}_{\,\dform{c}},\ \ 
\hat{\X}^\form{C}_{\,\dform{c}} = \X^\form{C}_{\,\dform{c}}
\end{align*}
and the normal tractor covariant derivative has the form
$$
\na_b
\begin{pmatrix} 
\si^\form{c} \\ \rh^{\dform{c}} 
\end{pmatrix} = 
\begin{pmatrix} 
\na_b \si^\form{c} + \rh^{\dform{c}} \de_b{}^{c^1} \\ 
\na_b \rh^{\dform{c}} - \ell P_{bc^1} \si^\form{c}
\end{pmatrix}, 
\ \ \text{i.e.} \ \
\na_b^{} \Y^\form{C}_{\,\form{c}} = 
-\ell\, P_{bc^1} \X^\form{C}_{\,\dform{c}},\ 
\na_b^{} \X^\form{C}_{\,\dform{c}} = \Y^\form{C}_{[b\dform{c}]}
$$

\begin{example} \label{skew}
We shall compute the sequence for the tractor bundle
$\cE^\form{C}$, $\form{C} = \form{C}^\ell$, i.e.\
$\cE^\form{C} \stackrel{\tilde{d}}{\to} \ldots 
\stackrel{\tilde{d}}{\to} \cE_{\form{a}^n}{}^\form{C}$.
Since the filtration of $\cE^\form{C}$ has level 2, it follows 
immediately from the construction of $\tilde{d}$ that
$(\tilde{d}F)_{a^0\form{a}} = (d^\na F)_{a^0\form{a}}{}^\form{C}
+ (\Box_{k+1})^{-1} (\pa^* d^\na d^\na F)_{a^0\form{a}}{}^\form{C}$
for every $F_\form{a}{}^\form{C} \in \cE_\form{a}{}^\form{C}$.
(In particular, the difference between $d^\na$ and $\tilde{d}$ is algebraic in 
this case.) 

Let us compute $\tilde{d}$ in details. Assume
$F_\form{a}{}^\form{C} = \Y^\form{C}_{\,\form{c}} \si_\form{a}{}^{\form{c}} + 
\X^\form{C}_{\,\dform{c}} \rh_\form{a}{}^{\dform{c}}$. Then
\begin{align*}
(d^\na d^\na F)_{a^{-1}a^0\form{a}}{}^\form{C} &= 
\frac{1}{2} \Om_{a^{-1}a^0} \sharp F_\form{a}{}^\form{C} =
\frac{1}{2}\ell \Om_{a^{-1}a^0}{}^{C^1}{}_P F_\form{a}{}^{P\dform{C}} =\\
&=\frac{1}{2} \ell
\Y^\form{C}_{\,\form{c}} W_{a^{-1}a^0}{}^c{}_p \si_\form{a}{}^{p\dform{c}}
+ \X^\form{C}_{\,\dform{c}} \bar{\rh}_{a^{-1}a^0\form{a}}{}^{\dform{c}}
\end{align*}
for some section $\bar{\rh}$ which we shall not need explicitly. 
Therefore 
\begin{align*}
(\pa^* d^\na d^\na F)_{a^0\form{a}}{}^\form{C} 
=& \frac{\ell^2}{2} X^{C^1} X_Q^r 
\Om_{[ra^0}{}^{[Q}{}_{|P|} F_{\form{a}]}{}^{|P|\dform{C}]}
=\frac{\ell^2}{2} \X^\form{C}_{\,\dform{c}} 
W_{[ra^0}{}^{[r}{}_{|p|} \si_{\form{a}]}{}^{|p|\dform{c}]} = \\
=& \frac{\ell}{2(k+2)} \X^\form{C}_{\,\dform{c}} \bigl[
-(\ell-1) W_{pr}{}^{c^2}{}_{a^0} 
\si_\form{a}{}^{pr\ddform{c}} + 
k W_{a^0a^1}{}^r{}_p \si_{r\dform{a}}{}^{p\dform{c}} \bigr].
\end{align*}

It remains to apply $(\Box_{k+1})^{-1}$. Note the map
$\pa^* d^\na d^\na: \cE_\form{a}{}^\form{C} \to \cE_{a^0\form{a}}{}^\form{C}$
has values in the (completely reducible) subbundle
$\cE_{a^0\form{a}}{}^{\dform{c}}(-\ell) \subseteq \cE_{a^0\form{a}}{}^\form{C}$,
cf.\ the precious display. 
Irreducible components of this subbundle are bundles
$\text{tf}[\cE_{\form{b}^{k+2-i}}{}^{\form{d}^{\ell-i}}](-\ell)$, 
$1 \leq i \leq \min \{\ell,k+2\}$ where 
the notation $\text{tf}[..]$ denotes the trace--free part of the enclosed 
bundle. The Laplace-Kostant operator $\Box_{k+1}$ on 
$\text{tf}[\cE_{\form{b}^s}{}^{\form{d}^t}](-\ell)$ acts by
$A_s^t(\ell) := \frac{1}{s+1}[n-s-t+1 + (l-t)(n-s)]$. Note the computation
is rather simple if we consider
$\text{tf}[\cE_{\form{b}^s}{}^{\form{d}^t}](-\ell)$ as the irreducible
invariant subbundle of $\cE^{\form{D}^t(E_1 \ldots E_{l-t})}$
and then follow \ref{standard}.
Also note $A_s^t(\ell)$ is always nonzero. This of course follows by general 
means but can be verified directly since 
$\text{tf}[\cE_{\form{b}^s}{}^{\form{d}^t}] \not= \{0\}$ if and only if 
$s+t \leq n$.

\begin{proposition}
The operator 
$\tilde{d}: \cE_\form{a}{}^\form{C} \to \cE_{a^0\form{a}}{}^\form{C}$
in the projective geometry has the form 
$$
(\tilde{d} F)_{a^0\form{a}}{}^\form{C} = (d^\na F)_{a^0\form{a}}{}^\form{C}
- \frac{\ell^2}{2} \sum_{i=1}^{\min\{\ell,k+2\}}
\frac{1}{A_{k+2-i}^{\ell-i}(\ell)} 
\Proj_{k+2-i}^{\ell-i} \X^\form{C}_{\,\dform{c}} 
W_{[ra^0}{}^{[r}{}_{|p|} \si_{\form{a}]}{}^{|p|\dform{c}]}
$$
where $\si_{\form{a}}{}^{\form{c}} = 
\X^{\,\form{c}}_{\form{C}} F_{\form{a}}{}^\form{C}$, 
$\X^{\,\form{c}}_{\form{C}} = X^{\,c^1}_{C^1} \ldots X^{\,c^\ell}_{C^\ell}$
and $\Proj_s^t: \cE_{\form{a}^{s+i}}{}^{\form{c}^{t+i}}(\ell) 
\to \text{tf}[\cE_{\form{a}^s}{}^{\form{c}^t}](\ell)$, $i \geq 0$ is the 
projection.
\qed
\end{proposition}

The operator $\tilde{d}$ simplifies in special cases $\ell=1$ and $k=0$. 
First assume $\ell=1$. Then
$(\pa^* d^\na d^\na)_{a^0\form{a}}{}^C = 
\frac{k}{2(k+2)} X^\form{C}
W_{a^0a^1}{}^r{}_p \si_{r\dform{a}}^{}{}^{p}$ has values in the irreducible
subbundle $\cE_{a^0\form{a}}(-\ell)$ of $\cE_{a^0\form{a}}{}^C$. We computed
$\Box_{k+1}$ acts by $\frac{n-(k+1)}{k+2}$ on this subbundle. Inverting this
scalar, we obtain the result
$$
(\tilde{d}F)_{a^0\form{a}}{}^C = (d^\na F)_{a^0\form{a}}{}C
+ \frac{k}{2(n-k-1)} X^C
W_{a^0a^1}{}^r{}_p \si_{r\dform{a}}{}^{p}.
$$

Now assume $k=0$. Then 
$(\pa^* d^\na d^\na F)_a{}^\form{C} 
= -\frac{\ell(\ell-1)}{4} \X^\form{C}_{\,\dform{c}} 
W_{pr}{}^{c^2}{}_a \si_{}^{pr\ddform{c}}$
has values in the trace--free (thus irreducible) part of the
subbundle $\cE_a{}^{\dform{c}}(-\ell)$.
Since $\Box_{k+1}$ acts on the trace--free part of 
$\cE_a{}^{\dform{c}}(-\ell) \subseteq \cE_a{}^\form{C}$ by 
$\frac{n-\ell}{2}$, the resulting formula is 
$$
(\tilde{d}F)_a{}^\form{C} = (d^\na F)_a{}^\form{C}
+ \frac{\ell(\ell-1)}{2(n-\ell)} \X^\form{C} 
W_{pr}{}^{c^2}{}_a \si^{pr\ddform{c}}.
$$
We claim $\tilde{d}$ actually coicides with the prolongation covariant 
derivative $\tilde{\na}$. To verify this, first observe 
$((\tilde{\na}-\na)F)_a{}^\form{C} \in \Im \pa^*$ by the constrution of 
$\tilde{d} = \tilde{\na}$. Thus it remains to verify 
$(d^{\tilde{\na}} \tilde{\na} F)_{a^{-1}a^0}{}^\form{C} \in \Ker \pa^*$.
But since $(d^{\na} \tilde{\na} F)_{a^{-1}a^0}{}^\form{C} \in \Ker \pa^*$
(again by the constrution of $\tilde{d} = \tilde{\na}$) and 
$d^{\tilde{\na}} - d^\na: \cE_{a^0} \to \ker \pa^* \subseteq
\cE_{a^{-1}a^0}{}^\form{C}$, cf.\ the last term in the previous display,
the claim follows. Using the matrix notation, $\tilde{\na} = \tilde{d}$
has the form
$$
\tilde{\na}_a
\begin{pmatrix}
\si^\form{c} \\ \rh^{\dform{c}}
\end{pmatrix}
= \na_a
\begin{pmatrix}
\si^\form{c} \\ \rh^{\dform{c}}
\end{pmatrix}
+ \frac{\ell(\ell-1)}{2(n-\ell)}  
\begin{pmatrix}
0 \\ W_{pr}{}^{c^2}{}_a \si^{pr\ddform{c}}
\end{pmatrix}. 
$$

Finally note $\cE^{\form{C}} \cong \cE_{\form{D}}$ (using the 
tractor volume form) for $\form{C} = \form{C}^\ell$ and 
$\form{D} = \form{D}^{n-\ell+1}$. The case $\ell = n-1$ 
(i.e.\ $\form{D} = \form{D}^2$) was solved in \cite{EaNotes},
where the prolongation of the corresponding BGG operator 
$\cE_a(2) \to \cE_{(ab)}$ (explicitly $f_a \mapsto \na_{(a} f_{b)}$)
is constructed. They construct the prolongation as the tractor covariant 
derivative $D_a: \cE_{\form{D}^2} \to \cE_{a\form{D}^2}$, cf.\
\cite{hsss}. Since $D_a-\na_a: \cE_{\form{D}^2} \to \im \pa^*$
(this follows from the formula for $D_a$ in p.\ 9, \cite{EaNotes} after 
a short computation) and the curvature of 
$(D_aD_b - D_bD_a): \cE_{\form{D}^2} \to \Ker \pa^*$ (this is obvious
form the formula for $D_aD_b - D_bD_a$ on $\cE_{\form{D}^2}$ on the same page)
we conclude $D_a = \tilde{\na}_a$, cf.\ \ref{thm-prolongcon}.
\end{example}

\begin{example} \label{symm}
Here we discuss the bundle 
$\cE^{(AB)} = \cE^{(ab)}(-2) \lpl \cE^a(-2) \lpl \cE(-2)$. 
Consider a section $F_\form{a}{}^{BC} \in \cE_\form{a}{}^{(BC)}$, expanded in
the basis of injectors as
$F_\form{a}{}^{BC} = Y^{(B}_{\,\ b} Y^{C)}_{\, c} \si_\form{a}{}^{bc} 
+ X^{(B} Y^{C)}_{\, c} \rh_\form{a}{}^c + X^B X^C \nu_\form{a}$. Then 
\begin{align*}
&(d^\na d^\na  F)_{a^{-1}a^0\form{a}}{}^{BC} = 
\frac{1}{2} \Om_{a^{-1}a^0} \sharp F_\form{a}{}^{BC} =
\Om_{a^{-1}a^0}{}^{(B}{}_P F_\form{a}{}^{C)P} =\\
&= Y^{(B}_{\,\ b} Y^{C)}_{\, c} W_{a^{-1}a^0}{}^{(b}{}_p \si_\form{a}{}^{c)p}
+X^{(B} Y^{C)}_{\, c} \bigl[
\frac{1}{2} W_{a^{-1}a^0}{}^c{}_p \rh_\form{a}{}^p 
- A_{a^{-1}a^0p} \si_\form{a}{}^{cp} \bigr]
+ X^B X^C \bar{\nu}_\form{a}
\end{align*}
for some section $\bar{\nu}$. Applying $\pa^*$ we obtain
\begin{align*}
(\pa^* d^\na d^\na F)_{a^0\form{a}}{}^{BC}
=& 2 X^{(B} Y^{C)}_{\, c} 
W_{[ra^0}{}^{(r}{}_{|p|} \si_{\form{a}]}{}^{c)p} \\
&+ X^B X^C \bigl[
\frac{1}{2} W_{[ra^0}{}^{r}{}_{|p|} \rh_{\form{a}]}{}^{p} 
- A_{[ra^0|p|} \si_{\form{a}]}{}^{pr} 
\bigr].
\end{align*}
The filtration degree of $\cE^{(AB)}$ is 3 and so
the construction of $\tilde{d}$ will require (at most)
2 steps. In the first step we put
$d' := d^\na + (\Box_{k+1}^{XY})^{-1} \pa^* d^\na d^\na: 
\cE_\form{a}{}^{BC} \to \cE_{a^0\form{a}}{}^{BC}$ where
$\Box_{k+1}^{XY}$ denotes $\Box_{k+1}$ restricted to the subquotient
$\cE_\form{a}{}^{c}(-2)$ of $\cE_\form{a}{}^{(BC)}$ which corresponds to the 
injector 
$X^{(B} Y^{C)}_{\, c}: \cE_\form{a}{}^{c}(-2) \hookrightarrow 
\cE_\form{a}{}^{(BC)}$. Note this subquotient has two irreducible components
but we need only the trace--free part since 
$W_{[ra^0}{}^{(r}{}_{|p|} \si_{\form{a}]}{}^{c)p}$ is trace--free. 
A short computation reveals 
$\pa^*\pa = \Box_1$ acts on the corresponding subquotient
of $\cE_a{}^{(BC)}$ by $\frac{n-k}{k+2}$. Hence 
\begin{align} \label{d'}
\begin{split} 
(d'F)_{a^0\form{a}}{}^{BC} = \na_{a^0} &F_\form{a}{}^{BC} 
- \frac{k+2}{n-k} \bigl[ 2 X^{(B} Y^{C)}_{\, c} 
W_{[ra^0}{}^{(r}{}_{|p|} \si_{\form{a}]}{}^{c)p} \\
&\qquad + X^B X^C \bigl(
\frac{1}{2} W_{[ra^0}{}^{r}{}_{|p|} \rh_{\form{a}]}{}^{p} 
- A_{[ra^0|p|} \si_{\form{a}]}{}^{pr} \bigr)
\bigr].
\end{split}
\end{align}
Further computation reveals
\begin{align*} 
\begin{split} 
&(d^\na d'F)_{a^{-1}a^0\form{a}}{}^{BC} = 
(d^\na d^\na F)_{a^{-1}a^0\form{a}}{}^{BC} 
- \frac{k+2}{n-k} \Bigl[ 2 Y^{(B}_{\,a^{-1}} Y^{C)}_{\, c} 
W_{[ra^0}{}^{(r}{}_{|p|} \si_{\form{a}]}{}^{c)p} \\
&\ + 2 X^{(B} Y^{C)}_{\,c} \bigl(
+\frac{1}{2} \de_{a^{-1}}^c W_{[ra^0}{}^{r}{}_{|p|} \rh_{\form{a}]}{}^{p} 
- \de_{a^{-1}}^c A_{[ra^0|p|} \si_{\form{a}]}{}^{pr} \\
&\qquad\qquad\qquad\! + 
\na_{a^{-1}} W_{[ra^0}{}^{(r}{}_{|p|} \si_{\form{a}]}{}^{c)p}  
\bigr) \Bigr] + X^BX^C \ga_{a^{-1}a^0\form{a}}.
\end{split}
\end{align*}
for some section $\ga_{a^{-1}a^0\form{a}} \in \cE_{a^{-1}a^0\form{a}}(-2)$ and 
\begin{align*}
\begin{split} 
(\pa^* d^\na d'F)_{a^0\form{a}}{}^{BC} = 
& - \frac{1}{n-k} X^{B} X^{C} \Bigl[
2 \na_s^{} W_{[ra^0}{}^{(r}{}_{|p|} \si_{\form{a}]}{}^{s)p} \\
&\qquad\quad  + (n-k-2) \bigl( 
\frac{1}{2} W_{[ra^0}{}^{r}{}_{|p|} \rh_{\form{a}]}{}^{p} 
- A_{[ra^0|p|} \si_{\form{a}]}{}^{pr} \bigr)
\Bigr].
\end{split}
\end{align*}
The previous displays shows that
$(\pa^*d^\na d' F)_{a^0\form{a}}{}^{BC}$ is the section of the
subbundle $\cE_{a^0\form{a}}(-2) \subseteq \cE_{a^0\form{a}}{}^{BC}$.
Since $\Box_{k+1}$ acts on this sunbundle by $\frac{2(n-k-1)}{k+2}$,
we obtain the result
$\tilde{d} := d' - \frac{k+2}{2(n-k-1)} \pa^* d^\na d'$.

\begin{proposition} \label{prBC}
The operator 
$\tilde{d}: \cE_\form{a}{}^\form{(BC)} \to \cE_{a^0\form{a}}{}^\form{(BC)}$
in the projective geometry has the form 
\begin{align*}
(\tilde{d} F)_{a^0\form{a}}{}^\form{BC} = &\na_{a^0} F_\form{a}{}^{BC} 
- \frac{k+2}{n-k} \Bigl[ 2 X^{(B} Y^{C)}_{\, c} 
W_{[ra^0}{}^{(r}{}_{|p|} \si_{\form{a}]}{}^{c)p} \\
& - \frac{1}{2(n-k-1)} X^B X^C \bigl[
2 \na_s^{} W_{[ra^0}{}^{(r}{}_{|p|} \si_{\form{a}]}{}^{s)p} \\
&\qquad\qquad\qquad  - (n-k) \bigl( 
\frac{1}{2} W_{[ra^0}{}^{r}{}_{|p|} \rh_{\form{a}]}{}^{p} 
- A_{[ra^0|p|} \si_{\form{a}]}{}^{pr} \bigr)
\bigr] \Bigr].
\end{align*}
where $\si_{\form{a}}{}^{bc} = 
X^{\,b}_{B} X^{\,c}_{C} F_{\form{a}}{}^{BC}$ and 
$\rh_{\form{a}}{}^{b} = 
2X^{\,b}_{B} Y_{C}^{} F_{\form{a}}{}^{BC}$.
\qed
\end{proposition}

We shall discuss the case $k=0$ in more details. Then
the formula in  Proposition \ref{prBC} simplifies to
\begin{align*}
(\tilde{d} F)_a{}^{BC} =& \na_a F^{BC} 
- \frac{2}{n} X^{(B} Y^{C)}_{\, c} W_{ra}{}^{c}{}_{p} \si^{rp}  \\
&+\frac{1}{n} X^B X^C \bigl( 2 A_{rap} \si^{pr} 
+ \frac{1}{n-1} W_{ra}{}^s{}_p \na_s \si^{rp} \bigr]. 
\end{align*}
This means $\tilde{d}$ is not a covariant derivative on $\cE^{(BC)}$ as the
term $W_{ra}{}^s{}_p \na_s \si^{rp}$ is not algebraic in 
$F^{BC}$, i.e.\ $\tilde{d} \not= \widetilde{\na}$ in this case. To compute 
$\widetilde{\na}$
explicitly, assume $k=0$ and put $\na' := d'$ (this is a covariant derivative
on $\cE^{(BC)}$). That is,
$\na'_a F^{BC} = \na_a F^{BC} -\frac{2}{n} (\Ps F )_a{}^{BC}$,
where the homomorphism $\Ps_a: \cE^{(BC)} \to \cE_a{}^{(BC)}$ is given
by the formula \nn{d'}, i.e.
$(\Ps F )_a{}^{BC} = 
X^{(B} Y^{C)}_{\, c} W_{ra}{}^{c}{}_{p} \si^{rp} 
- X^B X^C A_{rap} \si^{pr}.$
Extending $\Ps_{a^0}$ to an endomorphism 
$\cE_{{a^1}}{}^{(BC)} \to \cE_{a^0a^1}{}^{(BC)}$, an easy computation 
shows
$$
(\Ps \na'F)_{a^0a^1}{}^{BC} = 
X^{(B} Y^{C)}_{\, c} \bigl[ 
W_{ra^0}{}^{c}{}_{p} \na_{a^1} \si^{rp} 
- \frac{3}{2} W_{a^0a^1}{}^{c}{}_{p} \rh^p \bigr]
+ X^B X^C \bar{\nu} 
$$
for some $\bar{\nu} \in \cE(-2)$. 
Therefore 
$(\pa^* \Ps \na'F)_a^{BC} = 
- \frac{1}{2} X^B X^C W_{ra}{}^{c}{}_p \na_c\si^{rp}$
and we finally obtain
$(\pa^* d^{\na'} \na' F)_a{}^{BC} = 
(\pa^* d^{\na} \na' F)_a{}^{BC} - \frac{2}{n} (\pa^* \Ps \na' F )_a^{BC}=0$.
Since the left hand side is the curvature of $\na'$ (applied to $F^{BC}$),
this curvature is a map $\cE^{(BC)} \to \Ker \pa^*$.
Thus we verified $\widetilde{\na} = \na'$, cf.\ Theorem \ref{thm-prolongcon}.
Rewritting $\widetilde{\na}$ in the matrix notation, we obtain
$$
\widetilde{\na}_a
\begin{pmatrix}
\si^{bc} \\ \rh^c \\ \nu
\end{pmatrix}
= \na_a
\begin{pmatrix}
\si^{bc} \\ \rh^c \\ \nu
\end{pmatrix}
- \frac{2}{n} 
\begin{pmatrix}
0 \\ W_{ra}{}^{c}{}_p \si^{pr} \\ - A_{rap} \si^{pr}
\end{pmatrix}. 
$$

Note $\widetilde{\na}_a$ provides the prolongation of the corresponding 
(first order) BGG operator from $\cE^{(bc)_0}(-2)$ to the totally 
trace--free part of 
$\cE_a{}^{(bc)}(-2)$. The same problem was solved in \cite{EaMa}
in terms of the connection defined by (3.6) or the left hand side of 
(5.2) there. Let us denote this connection on $\cE^{(BC)}$ by $D_a$.
Note the formula for $D_a$ differs from $\widetilde{\na}_a$ in the middle
term of the last matrix in the previous display: this term is
$-\frac{2}{n} W_{ra}{}^{c}{}_p \si^{pr}$ for $\widetilde{\na}_a$ whereas
$\frac{1}{n} W_{ra}{}^{c}{}_p \si^{pr}$ in the case of $D_a$, cf.\
\cite[(3.6)]{EaMa}. The reason is purely notational, specifically in the 
choice of the projectors. If one replaces $X^{(B}Y^{C)}_{\,c}$ by
$-\frac{1}{2} X^{(B}Y^{C)}_{\,c}$ -- which means e.g.\
$F_\form{a}{}^{BC} = Y^{(B}_{\,\ b} Y^{C)}_{\, c} \si^{bc} 
+ (-\frac{1}{2} X^{(B} Y^{C)}_{\,c}) \rh^c + X^B X^C \nu$ --
both terms will coincide. Note also that formulas for $\na_a$ and the normal
covariant derivative defined in the display preceding to \cite[Theorem 5.1]{EaMa}
coincide after the change of projectors. This confirms the results here 
coincide with those in \cite{EaMa}.
\end{example}

\section{Conformal geometry}

\subsection{Conformal geometry and tractor calculus}

We summarise here some notation and background.  Further details may
be found in \cite{GoPetLap}.  Let $M$ be a smooth manifold of
dimension $n\geq 3$. Recall that a {\em conformal structure\/} of
signature $(p,q)$ on $M$ is a smooth ray subbundle $\cQ\subset
S^2T^*M$ whose fiber over $x$ consists of conformally related
signature-$(p,q)$ metrics at the point $x$. Sections of $\cQ$ are
metrics $g$ on $M$. So we may equivalently view the conformal
structure as the equivalence class $[g]$ of these conformally related
metrics.  The principal bundle $\pi:\cQ\to M$ has structure group
$\R_+$, and so each representation ${\R}_+ \ni x\mapsto x^{-w/2}\in
{\rm End}(\R)$ induces a natural line bundle on $ (M,[g])$ that we
term the conformal density bundle $E[w]$. We shall write $\cE[w]$ for
the space of sections of this bundle.  We write $\cE^a$ for the space
of sections of the tangent bundle $TM$ and $\cE_a$ for the space of
sections of $T^*M$. The indices here are abstract in the sense of
\cite{ot} and we follow the usual conventions from that source. So for
example $\cE_{ab}$ is the space of sections of $\otimes^2T^*M$.  Here
and throughout, sections, tensors, and functions are always smooth.
When no confusion is likely to arise, we will use the same notation
for a bundle and its section space.

We write $\bg$ for the {\em conformal metric}, that is the
tautological section of $S^2T^*M\otimes E[2]$ determined by the
conformal structure. This is used to identify $TM$ with
$T^*M[2]$.  For many calculations we employ abstract indices in an
obvious way.  Given a choice of metric $ g$ from $[g]$,
we write $ \nabla$ for the corresponding Levi-Civita connection. With
these conventions the Laplacian $ \Delta$ is given by
$\Delta=\bg^{ab}\na_a\na_b= \na^b\na_b\,$. Here we are raising indices
and contracting using the (inverse) conformal metric. Indices will be
raised and lowered in this way without further comment.  Note $E[w]$
is trivialised by a choice of metric $g$ from the conformal class, and
we also write $\na$ for the connection corresponding to this
trivialisation.  The coupled $\na_a$
preserves the conformal metric.

The curvature $R_{ab}{}^c{}_d$ of the Levi-Civita connection (the
Riemannian curvature) is given by $ [\na_a,\na_b]v^c=R_{ab}{}^c{}_dv^d
$ ($[\cdot,\cdot]$ indicates the commutator bracket).  This can be
decomposed into the totally trace-free Weyl curvature $C_{abcd}$ and a
remaining part described by the symmetric {\em Schouten tensor}
$\Rho_{ab}$, according to
\begin{equation}\label{csplit}
R_{abcd}=C_{abcd}+2\bg_{c[a}\Rho_{b]d}+2\bg_{d[b}\Rho_{a]c},
\end{equation}
 where
$[\cdots]$ indicates antisymmetrisation over the enclosed indices.
The Schouten tensor is a trace modification of the Ricci tensor
$\Ric_{ab}=R_{ca}{}^c{}_b$
and vice versa: $\Ric_{ab}=(n-2)\Rho_{ab}+\J\bg_{ab}$,
where we write $ \J$ for the trace $ \Rho_a{}^{a}$ of $ \Rho$.  The {\em
Cotton tensor} is defined by
$
A_{abc}:=2\nabla_{[a}\Rho_{b]c} .
$
Via the Bianchi identity this is related to the divergence of the Weyl
tensor as follows:
\begin{equation}\label{bi1} (n-3)A_{abc}=\nabla^d
C_{dcab} .
\end{equation}
Finally we put
\begin{equation} \label{Bach}
B_{ab} = \na^p A_{pab} + \Rho^{pq} C_{paqb} \in \cE_{(ab)_0}[-2].
\end{equation}
In the dimension $n=4$, this is the conformally invariant \idx{Bach} tensor.

Under a {\em conformal transformation} we replace a choice of metric
$g$ by the metric $\hat{g}=e^{2\Up} g$, where $\Up$
 is a smooth
function. We recall that, in particular, the Weyl curvature is
conformally invariant $\widehat{C}_{abcd}=C_{abcd}$.
With $\Up_a: = \na_a \Up$, the Schouten tensor transforms according to
\begin{equation}\label{Rhotrans}
\textstyle \widehat{\Rho}_{ab}=\Rho_{ab}-\na_a \Up_b +\Up_a\Up_b
-\frac{1}{2} \Up^c\Up_c \bg_{ab} .
\end{equation}

Explicit formula for the corresponding transformation of
the Levi-Civita connection and its curvatures are given in e.g.\
\cite{thomass,GoPetLap}. From these, one can easily compute the
transformation for a general valence (i.e.\ rank) $s$ section
$f_{bc \cdots d} \in \cE_{bc \cdots d}[w]$ using the Leibniz rule:
\begin{equation} \label{grad_trans_gen}
\begin{split}
  \hat{\na}_{\bar{a}} f_{bc \cdots d}
  =&  \na_{\bar{a}} f_{bc \cdots d} + (w-s) \Up_{\bar{a}} f_{bc \cdots d}
     -\Up_{b} f_{\bar{a}c \cdots d} \cdots
     -\Up_{d} f_{bc \cdots \bar{a}} \\
   & +\Up^p f_{p c \cdots d} \bg_{b\bar{a}} \cdots
     +\Up^p f_{bc \cdots p} \bg_{d\bar{a}}.
\end{split}
\end{equation}

We next define the standard tractor bundle over $(M,[g])$.
It is a vector bundle of rank $n+2$ defined, for each $g\in[g]$,
by  $[\cE^A]_g=\cE[1]\oplus\cE_a[1]\oplus\cE[-1]$.
If $\wh g=e^{2\Up}g$, we identify
 $(\alpha,\mu_a,\tau)\in[\cE^A]_g$ with
$(\wh\alpha,\wh\mu_a,\wh\tau)\in[\cE^A]_{\wh g}$
by the transformation
\begin{equation}\label{transf-tractor}
 \begin{pmatrix}
 \wh\alpha\\ \wh\mu_a\\ \wh\tau
 \end{pmatrix}=
 \begin{pmatrix}
 1 & 0& 0\\
 \Up_a&\delta_a{}^b&0\\
- \tfrac{1}{2}\Up_c\Up^c &-\Up^b& 1
 \end{pmatrix}
 \begin{pmatrix}
 \alpha\\ \mu_b\\ \tau
 \end{pmatrix} .
\end{equation}
It is straightforward to verify that these identifications are
consistent upon changing to a third metric from the conformal class,
and so taking the quotient by this equivalence relation defines the
{\em standard tractor bundle} $\cE^A$ over the conformal manifold.
On a
conformal structure of signature $(p,q)$, the bundle $\cE^A$ admits an
invariant metric $ h_{AB}$ of signature $(p+1,q+1)$ and an invariant
connection, which we shall also denote by $\nabla_a$, preserving
$h_{AB}$. Up to an isomorphism this the unique {\em normal conformal
tractor connection} and so induces normal
connection on $\bigotimes \cE^A$ that will be denoted $\na_a$
and termed the (normal) tractor connection.  In a conformal scale $g$,
the metric $h_{AB}$ and $\na_a$ on $\cE^A$ are given by
\begin{equation}\label{basictrf}
 h_{AB}=\begin{pmatrix}
 0 & 0& 1\\
 0&\bg_{ab}&0\\
1 & 0 & 0
 \end{pmatrix}
\text{ and }
\nabla_a\begin{pmatrix}
 \alpha\\ \mu_b\\ \tau
 \end{pmatrix}
 =
\begin{pmatrix}
 \nabla_a \alpha-\mu_a \\
 \nabla_a \mu_b+ \bg_{ab} \tau +\Rho_{ab}\alpha \\
 \nabla_a \tau - \Rho_{ab}\mu^b  \end{pmatrix}.
\end{equation}
It is readily verified that both of these are conformally well-defined,
i.e., independent of the choice of a metric $g\in [g]$.  Note that
$h_{AB}$ defines a section of $\cE_{AB}=\cE_A\otimes\cE_B$, where
$\cE_A$ is the dual bundle of $\cE^A$. Hence we may use $h_{AB}$ and
its inverse $h^{AB}$ to raise or lower indices of $\cE_A$, $\cE^A$ and
their tensor products.

In computations, it is often useful to introduce
the `projectors' from $\cE^A$ to
the components $\cE[1]$, $\cE_a[1]$ and $\cE[-1]$ which are determined
by a choice of scale.
They are respectively denoted by $X_A\in\cE_A[1]$,
$Z_{Aa}\in\cE_{Aa}[1]$ and $Y_A\in\cE_A[-1]$, where
 $\cE_{Aa}[w]=\cE_A\otimes\cE_a\otimes\cE[w]$, etc.
 Using the metrics $h_{AB}$ and $\bg_{ab}$ to raise indices,
we define $X^A, Z^{Aa}, Y^A$. Then we
 see that
$
Y_AX^A=1,\ \ Z_{Ab}Z^A{}_c=\bg_{bc} ,
$
and all other quadratic combinations that contract the tractor
index vanish.
In \eqref{transf-tractor} note that
$\wh{\alpha}=\alpha$ and hence $X^A$ is conformally invariant.
Reformulating \nn{basictrf}, we obtain
$$
\na_a^{} Y_B^{} = Z_B^b P_{ab}, \quad 
\na_a^{} Z_B^b = - Y_B^{} \de^b_a - X_B^{} P_a{}^b 
\quad \text{and} \quad
\na_a^{} X_B^{} = Z_B^b \bg_{ab}.
$$

Given a choice of $g\in [g]$,
the {\em tractor-$D$ operator}
$
D_A\colon\cE_{B \cdots E}[w]\to\cE_{AB\cdots E}[w-1]
$
is defined by
\begin{equation}\label{Dform}
D_A V:=(n+2w-2)w Y_A V+ (n+2w-2)Z_{Aa}\nabla^a V -X_A\Box V,
\end{equation}
 where $\Box V :=\Delta V+w \J V$.  This is
 conformally invariant, as can be checked directly using the formula
 above. 

The curvature $ \Omega$ of the tractor connection
is defined on $\cE^C$ by
$
[\na_a,\na_b] V^C= \Omega_{ab}{}^C{}_EV^E ~.
$
 Using
\eqref{basictrf} and the formulae for the Riemannian curvature yields
\begin{equation}\label{tractcurv}
\Omega_{abEF}= Z_E^{\,e} Z_F^{\,f} C_{abef}-2X_{[E}^{}Z_{F]}^{\,f} A_{abf}
\in \cE_{[ab][EF]} = \cE_{[ab]} \otimes \cA
\end{equation}
where $\cA = \cE_{[EF]}$ is the conformal adjoint tractor bundle.
We shall write $\Om_{ab} \sharp F_C$ or $(\Om \sharp F)_{abC}$ 
for the curvature action 
$(\na_a\na_b - \na_b\na_a) F_C = - \Om_{ab}{}^D{}_C F_D$.

\vspace{1ex}

Using the notation developed above, the inclusions $\io$ and $\bar{\io}$
defined in \ref{KoLap} have he form
$-2Y_{[E} Z_{F]a^0}: \cE_\form{a} \stackrel{\bar{\io}}{\to} 
\cE_{a^0\form{a}[EF]}$ and 
$-2X_{[E}^{} X_{F]}^{a^1}: 
\cE_\form{a} \stackrel{\io}{\to} \cE_{\dform{a}[EF]}$.
(The scalar $-2$ is used for the sake of compatibility of 
$\pa$ and $\na$, cf.\ \cite{CSbook}.) Thus
\begin{align*}
&\pa: \cE_\form{a} \otimes \cT \ni f_\form{a} \mapsto
-2Y_{[E} Z_{F]a^0} f_\form{a} \stackrel{\sharp}{\longrightarrow} 
\cE_{a^0\form{a}} \otimes \cT \quad \text{and} \\
&\pa^*: \cE_\form{a} \otimes \cT \ni f_\form{a} \mapsto
-2 X_{[E}^{} Z_{F]}^{a^1} f_\form{a} \stackrel{\sharp}{\longrightarrow} 
\cE_{\dform{a}} \otimes \cT
\end{align*}
and we can easily compute $\Box_k$ on $\cE_\form{a} \otimes \cT$ 
using the tensorial action $\sharp$.

\begin{example}
We shall compute $\tilde{d}$ on forms twisted by $\cE_C$. 
Let $\form{a} = \form{a}^k$ and consider
$F_{\form{a}C} = Y_C^{} \si_\form{a} + Z_C^{\,c} \mu_{c\form{a}}
+X_C^{} \nu_\form{a} \in \cE_{\form{a}C}$. Then 
\begin{align*}
(d^\na d^\na F)_{a^{-1}a^0\form{a}}{}_C &= 
\frac{1}{2} \Om_{a^{-1}a^0} \sharp F_\form{a}{}_C =
\frac{1}{2} \Om_{a^{-1}a^0}{}_{C}{}^P F_\form{a}{}_P =\\
&=\frac{1}{2} Z_C^{\,c} \bigl[ 
C_{a^{-1}a^0c}{}^p \mu_{\form{a}p} + A_{a^{-1}a^0c} \si_\form{a} \bigr]
- X_C A_{a^{-1}a^0}{}^p \mu_{\form{a}p}
\end{align*}
hence $(\pa^* d^\na d^\na F)_{a^0\form{a}}{}_C = 
- \frac{k}{2(k+2)} X_C \bigl[ C_{a^0a^1}{}^{rp} \mu_{r\dform{a}p} 
+ A_{a^0a^1}{}^r \si_{r\dform{a}}\bigr]$.
This is a section of the subbundle 
$\cE_{a^0\form{a}}[-1] \subseteq \cE_{a^0\form{a}C}$ and one easily computes
$\Box_k$ acts on this (irreducible) subbundle by 
$-\frac{n-k-1}{k+2}$. Therefore 
$(\tilde{d}F)_{a^0\form{a}C} = \na_{a^0} F_{\form{a}C}
- \frac{k}{2(n-k-1)} X_C \bigl[ C_{a^0a^1}{}^{rp} \mu_{r\dform{a}p} 
+ A_{a^0a^1}{}^r \si_{r\dform{a}}\bigr]$ for $0 \leq k \leq n-1$ and
$\tilde{d} = d^\na$ for $k \geq n-1$. Finally note that 
the prolongation covariant derivative coincides with the normal one
for $k=0$, i.e.\ $\tilde{\na} = \na$ on $\cE_C$. 
\end{example}

\begin{example}
The computation of the prolongation covariant derivative is getting rather 
technical for more complicated bundles. We shall demonstrate it on
the prolongation covariant derivative $\wt{\na}$
on $\cE_{(BC)_0}$. (Note $\cE_{(BC)_0}$ and $\cE^{(BC)_0}$ are isomorphic 
using the tractor metric.) The computation consists of three steps:
we start with $\na$ and then define covariant derivatives $\ol{\na}$,
$\ol{\ol{\na}}$ and $\wt{\na}$. Taking a section 
$F_{BC} = Y_{(B} Y_{C)} \si + Y_{(B}^{} Z_{C)}^{\,c} \rh_c + 
Z_{(B}^{\,b} Z_{C)}^{\,c} \om_{bc} + X_{(B} Y_{C)} \nu + 
X_{(B}^{} Z_{C)}^{\,c} \mu_c + X_{(B} X_{C)} \ka$ we get
\begin{align*}
(d^\na & d^\na  F)_{a^0a^1BC} = 
\frac{1}{2} \Om_{a^0a^1} \sharp F_{BC} =
\frac{1}{2} \Om'_{a^0a^1BC}{}^{PQ} F_{PQ}^{} =\\
=& Y_{(B}^{} Z_{C)}^{\, c} \bigl[
\frac{1}{2} C_{a^0a^1c}{}^{p} \rh_p + A_{a^0a^1c} \si \bigr]
+ Z_{(B}^{\,b} Z_{C)}^{\, c} \bigl[
C_{a^0a^1(b}{}^{p} \om_{c)p} + \frac{1}{2} A_{a^0a^1(b} \rh_{c)} \bigr] \\
& - \frac{1}{2} X_{(B} Y_{C)} A_{a^0a^1}{}^p \rh_p
+ X_{(B}^{} Z_{C)}^{\, c} \bigl[
\frac{1}{2} C_{a^0a^1c}{}^{p} \mu_p - A_{a^0a^1}{}^p \om_{cp} 
+ \frac{1}{2} A_{a^0a^1c} \nu \bigr] \\
& - \frac{1}{2} X_B X_C A_{a^0a^1}{}^p \mu_p.
\end{align*}
where $\Om'_{a^0a^1BC}{}^{PQ} := 2\Om_{a^0a^1(B}{}^{(P} h_{C)}{}^{Q)}$.
Applying $\pa^*$ to the previous display we obtain 
$(\pa^* d^\na  d^\na F)_{a^1BC} 
= -2 \X_{(B}{}^{Pr} \Om_{|ra^1|C)}{}^Q F_{PQ}$ because 
$\Om_{a^0a^1EF}$ is $\pa^*$-closed (i.e.\ $\X_{A^0}{}^{Pp} \Om_{pa^1PA^1}=0$). 
We put ${{\Ps}}_{a^1BC}{}^{PQ} := -2 \X_{(B}{}^{Pr} \Om_{|ra^1|C)}{}^Q$.
Equivalently, ${{\Ps}}_{a^1BC}{}^{PQ}$ can be obtained by applying 
$\pa^*$ to the $\cE_{BC}$-factor of
$\Om'_{a^0a^1(BC)}{}^{PQ}$. This is exactly
the operator $\pa^*_V$ from \cite{hsss} since the notation therein means 
$V = \cE_{(BC)_0}$, $V^* = \cE^{(PQ)_0}$ and therefore $\Om'_{a^0a^1BC}{}^{PQ}
\in \cE_{a^0a^1} \otimes \End(V)$ is the curvature tensor of $\na_a$ on 
$V = \cE_{(BC)_0}$. We shall denote the operator $\pa^*_V$ by 
$\pa^*_{BC}: \cE_{a^0a^1BC}{}^{PQ} \to \cE_{a^1BC}{}^{PQ}$ here. 
Thus we have ${{\Ps}}_{a^1BC}{}^{PQ} = 
\frac{1}{2} (\pa^*_{BC} \Om')_{a^1BC}{}^{PQ}$, explicitly
\begin{align} \label{Ps}
\begin{split}
{\Ps}_{a^1BC}{}^{PQ} =
& - Z_{(B}^{\,b} Z_{C)}^{\, c} \bigl[
X^{(P} Z^{Q)q} C_{a^1(bc)q} + X^P X^Q A_{a^1(bc)} \bigr] \\
& + X_{(B} Z_{C)}^{\, c} \bigl[
Z^{p(P} Z^{Q)q} C_{a^1pcq} + 2 X^{(P} Z^{Q)q} A_{a^1(cq)} \bigr] \\
& + X_{(B} X_{C)} Z^{p(P} Z^{Q)q} A_{pa^1q}.
\end{split}
\end{align}

Since $\frac{1}{2} C_{a^1(bc)}{}^{p} \rh_{p} + A_{a^1(bc)} \si$ is a section
of the Cartan component of the subquotient $\cE_{[a^1b]} \otimes \cE_c$ of
$\cE_{a^1(BC)_0}$ and 
$\Box_1$ acts on this subquotient by $-\frac{3}{2}$, we put
$\ol{\na}_a F_{BC} = \na_a F_{BC} + \frac{2}{3} {{\Ps}}_{aBC}{}^{PQ} F_{PQ}$
as the first ``approximation'' of $\wt{\na}$. We need to know 
$\na_{a^0} {\Ps}_{a^1BC}{}^{PQ}$ to compute the curvature
$\ol{\Om}_{a^0a^1BC}{}^{PQ}$ of $\ol{\na}$.
First, it easily follows from 
${{\Ps}}_{a^1BC}{}^{PQ} := -2 \X_{(B}{}^{Pr} \Om_{|ra^1|C)}{}^Q$ that
\begin{align*}
(d^\na & \Ps)_{a^0a^1BC}{}^{PQ} = \na_{a^0} {\Ps}_{a^1BC}{}^{PQ} =
-2 \na_{a^0}  \X_{(B}{}^{Pr} \Om_{|ra^1|C)}{}^Q = \\
&= - 2\Z_{(B}^{\ e^0 Pe^1} \bg_{|a^0e^0} \Om_{e^1a^1|C)}{}^Q
+ 2\W_{(B}{}^P \Om_{|a^0a^1|C)}{}^Q
- \X_{(B}{}^{Pr} \na_{|r} \Om_{a^0a^1|C)}{}^Q
\end{align*}
since $\na_{a^{-1}} \Om_{a^0a^1CQ}=0$. 
Expanding the expressions in the previous display we obtain
\begin{align*}
&(d^\na {\Ps})_{a^0a^1BC}{}^{PQ} =
-\frac{3}{2} Y_{(B}^{} Z_{C)}^{\, c} \bigl[
X^{(P} Z^{Q)q} C_{a^0a^1cq} + X^P X^{Q} A_{a^0a^1c} \bigr] \\
& +\frac{3}{2} X_{(B} Y_{C)} X^{(P} Z^{Q)q} A_{a^0a^1q} \\
& + Z_{(B}^{\ b} Z_{C)}^{\, c} \Bigl[
-2 Z^{p(P} Z^{Q)q} \bg_{a^0[b} C_{p]a^1cq} 
+\frac{1}{2} X^PX^Q \bigl( \na_b A_{a^0a^1c} + P_{b}{}^r C_{a^0a^1rc} \bigr)\\
&\qquad\qquad\quad  + X^{(P} Z^{Q)q} \bigl( -2 \bg_{a^0[b} A_{q]a^1c}
+ \frac{1}{2} \na_b C_{a^0a^1cq} - \bg_{b[c} A_{|a^0a^1|q]} \bigr) \Bigr] \\
& + X_{(B}^{} Z_{C)}^{\, c} \Bigl[ \;\frac{3}{2} Y^{(P} Z^{Q)q} C_{a^0a^1cq}
- X^{(P} Z^{Q)q} \bigl( \na_{(c} A_{|a^0a^1|q)} + P_{(c}{}^s C_{|a^0a^1s|q)}
\bigr) \\
&\qquad\qquad\quad +Z^{p(P} Z^{Q)q} \bigl( 2\bg_{a^0[c} A_{p]a^1q} 
- \frac{1}{2} \na_p C_{a^0a^1cq} + \bg_{p[c} A_{|a^0a^1|q]} \bigr) \Bigr] \\
& + X_B X_C \Bigl[  -\frac{3}{2} Y^{(P} Z^{Q)q} A_{a^0a^1q}
+ \frac{1}{2} Z^{p(P} Z^{Q)q} \bigl(
\na_p A_{a^0a^1q} + P_p{}^s C_{a^0a^1sq} \bigl) \Bigr]
\end{align*}
after some computation which uses the differential Bianchi identity, 
in particular the relation \cite[(29)]{GoSi}.
Now we need to apply $\pa^*_{BC}$ to the previous display. This yields
\begin{align*}
(\pa^*_{BC} \,& d^\na {\Ps})_{a^1BC}{}^{PQ} = 
\frac{3}{2} Z_{(B}^{\ b} Z_{C)}^{\, c} \bigl[ 
X^{(P} Z^{Q)q} C_{a^1(bc)q} + X^P X^Q A_{a^1(bc)} \bigr] \\
& + X_{(B}  Z_{C)}^{\, c} \bigr[ 
\frac{1}{2}(n-1) Z^{p(P} Z^{Q)q} C_{a^1(pq)c} - \frac{1}{2} X^PX^Q B_{a^1c} \\
&\qquad\qquad\quad\
+ X^{(P} Z^{Q)q} \bigl( (n-4) A_{q(a^1c)} - 3 A_{a^1(qc)} \bigr) \bigr] \\
& + X_BX_C \bigl[ \frac{1}{2} (n-1) Z^{p(P} Z^{Q)q} A_{a^1(pq)}
+ \frac{1}{2} X^{(P} Z^{Q)q} B_{a^1q} \bigr].
\end{align*}
We need to compute 
${\ol{\Ps}}_{a_1BC}{}^{PQ} = 
\frac{1}{2} (\pa^*_{BC} \ol{\Om})_{a^1BC}{}^{PQ}$ satisfying
${\ol{\Ps}}_{a^1BC}{}^{PQ} F_{PQ} =  
(\pa^* d^{\ol{\na}} \ol{\na} F)_{a^1BC}$.
Since $\ol{\na}_a F_{BC} = \na_a + \frac{2}{3} {{\Ps}}_{aBC}{}^{PQ}$
we have
$$
\frac{1}{2} \ol{\Om}_{a^0a^1BC}{}^{PQ} = 
\frac{1}{2} \Om'_{a^0a^1BC}{}^{PQ} 
+ \frac{2}{3} (d^\na \Ps)_{a^0a^1BC}{}^{PQ} 
+ \frac{4}{9} (\Ps \wedge \Ps)_{a^0a^1BC}{}^{PQ}
$$
where $(\Ps \wedge \Ps)_{a^0a^1BC}{}^{PQ} = 
{\Ps}_{a^0BC}{}^{RS} {\Ps}_{a^1RS}{}^{PQ}$. 
Since $\frac{1}{2} (\pa^*_{BC} \Om')_{a^1BC}{}^{PQ} = {\Ps}_{a^1BC}{}^{PQ}$ 
by definition of ${\Ps}$, applying $\pa^*_{BC}$ to the previous display yields
\begin{align} \label{Ps-}
\begin{split}
&\ol{{\Ps}}_{a_1BC}{}^{PQ} = 
\frac{1}{2} (\pa^*_{BC} \ol{\Om})_{a^1BC}{}^{PQ} =\\
&= {\Ps}_{a^1BC}{}^{PQ} 
+ \frac{2}{3} (\pa^*_{BC} d^\na \Ps)_{a^1BC}{}^{PQ}
+ \frac{4}{9} (\pa^* (\Ps \wedge \Ps))_{a^1BC}{}^{PQ} = \\
&= \frac{1}{3} X_{(B}  Z_{C)}^{\, c} \bigl[ 
(n-4) Z^{p(P} Z^{Q)q} C_{a^1(pq)c} + 2(n-4) X^{(P} Z^{Q)q} A_{q(a^1c)}
- X^PX^Q B_{a^1c} \bigr] \\
&\quad + \frac{1}{3} X_BX_C \bigl[ (n-4) Z^{p(P} Z^{Q)q} A_{a^1(pq)}
+ X^{(P} Z^{Q)q} B_{a^1q} \bigr]
+ \frac{4}{9} (\pa^*_{BC} (\Ps \wedge \Ps))_{a^1BC}{}^{PQ}.
\end{split}
\end{align}
where
\begin{equation} \label{*PsPs}
(\pa^*_{BC} (\Ps \wedge \Ps))_{a^1BC}{}^{PQ} = 
\frac{1}{2} X_{B}  X_C \bigl[ 
X^{(P} Z^{Q)q} C_{a^1}{}^{(rs)p} C_{qrsp} 
+ X^PX^Q C_{a^1}{}^{(rs)q} A_{qrs} \bigr].
\end{equation}

\begin{remark}
The section
$(\pa^* d^\na {\Ps})_{aBC}{}^{PQ}$ is closely related to the conformally
invariant curvature quantity
\begin{align*}
W_{\form{EF}} =
& (n-4) \Z_\form{E}^{\,\form{e}} \Z_\form{F}^{\,\form{f}} C_{\form{ab}}
- 2(n-4) \Z_\form{E}^{\,\form{e}} \X_\form{F}^{\,f} A_{\form{e}f} \\
& - 2(n-4) \X_\form{E}^{\,e} \Z_\form{F}^{\,\form{f}} A_{\form{f}e}
+ 4  \X_\form{E}^{\,e} \X_\form{F}^{\,f} B_{ef},
\end{align*}
cf.\ \cite{GoAdv} where all the form indices $\form{E}$, $\form{F}$, 
$\form{e}$, $\form{f}$ have the valence 2. In fact, one easily computes 
$(\pa^* d^\na {\Ps})_{aBC}{}^{PQ} = 
-\frac{1}{3} Z^R_{a} X_{(B} W_{C)}{}^{(P}{}_{R}{}^{Q)}$.
Since $(\pa^* d^\na {\Ps})_{aBC}{}^{PQ}$ coincides with 
$\ol{{\Ps}}_{aBC}{}^{PQ}$ up to the terms involving 
$C_{a^1}{}^{(rs)p} C_{qrsp}$ and $C_{a^1}{}^{(rs)q} A_{prs}$,
cf.\ \nn{Ps-}, conformal invariance of $W_{\form{EF}}$
verifies the invariance of the previous computations.
\end{remark}

Looking at the form of ${\ol{\Ps}}_{a_1BC}{}^{PQ} F_{PQ}$, we see
that we need the action of $\Box_1$ on the subquotient 
$\cE_{(a^1c)_0}$ of $\cE_{a_1BC}$ (corresponding to the injector
$X_{(B}^{} Z_{C)}^{\,c}$). A short computation reveals this is $-\frac{n}{2}$
hence the next ``approximation'' of $\wt{\na}$ will be the covariant 
derivative
$$
\ol{\ol{\na}}_a := \ol{\na}_a + \frac{2}{n} {\ol{\Ps}}_{aBC}{}^{PQ} =
\na_a + \frac{2}{3} {\Ps}_{aBC}{}^{PQ} + 
\frac{2}{n} {\ol{\Ps}}_{aBC}{}^{PQ}: \cE_{(PQ)} \to \cE_{a(BC)}.
$$
Now we need the curvature $\ol{\ol{\Om}}_{a^0a^1BC}{}^{PQ}$ of 
$\ol{\ol{\na}}_a$ and then to apply $\pa^*_{BC}$ on 
$\frac{1}{2} \ol{\ol{\Om}}_{a^0a^1BC}{}^{PQ}$. It follows from the 
definition of $\ol{\ol{\na}}_a$ that
\begin{align} \label{Om2}
\frac{1}{2} \ol{\ol{\Om}}_{a^0a^1BC}{}^{PQ} = 
& \frac{1}{2} \ol{\Om}_{a^0a^1BC}{}^{PQ}
+ \frac{2}{n} \na_{a^0} \ol{\Ps}_{a^1BC}{}^{PQ} 
+ \frac{4}{3n} \ol{\Ps}_{a^0BC}{}^{RS} {\Ps}_{a^1RS}{}^{PQ}
\end{align}
since $\ol{\Ps}_{a^0BC}{}^{RS} \ol{\Ps}_{a^1RS}{}^{PQ} = 
\Ps_{a^0BC}{}^{RS} \ol{\Ps}_{a^1RS}{}^{PQ}=0$.

The next step is to compute 
$\ol{\ol{\Ps}}_{a^1BC}{}^{PQ} := 
\frac{1}{2} (\pa^*_{BC} \ol{\ol{\Om}})_{a^1BC}{}^{PQ}$.
We apply $\pa^*_{BC}$ to the three terms on the right hand side of
\nn{Om2}. Firstly recall 
$\frac{1}{2} (\pa^*_{BC} \ol{\Om})_{a^1BC}{}^{PQ} 
= \ol{\Ps}_{a^1BC}{}^{PQ}$ by definition. Secondly, one gets
\begin{align*}
&(d^\na \ol{\Ps})_{a^0a^1BC}{}^{PQ} =  
\frac{1}{3} Z_{(B}^{\ b} Z_{C)}^{\, c} \Bigl[
(n-4) Z^{p(P} Z^{Q)q} \bg_{ba^0} C_{a^1(pq)c} \\
&\qquad\qquad\qquad\qquad + 2(n-4) X^{(P} Z^{Q)q} \bg_{ba^0} A_{q(a^1c)}  
- X^PX^Q \bg_{ba^0} B_{a^1c} \Bigr] \\
& + \frac{1}{3} X_{(B}^{} Z_{C)}^{\, c} \Bigl[ 
\frac{3}{2} (n-4) Y^{(P} Z^{Q)q} C_{a^0a^1qc}
-\frac{3}{2} (n-4) X^{(P} Y^{Q)} A_{a^0a^1c} \\
&\qquad + (n-4) Z^{p(P} Z^{Q)q} \bigl( 
\na_{a^0} C_{a^1(pq)c} +  2 \bg_{a^0(p} A_{q)(a^1c)} 
+ 2 \bg_{a^0c} A_{a^1(pq)} \bigr) \\
&\qquad + 2 X^{(P} Z^{Q)q} \bigl( (n-4) \na_{a^0} A_{q(a^1c)}
-(n-4) P_{a^0}{}^p C_{a^1(pq)c} + 2\bg_{a^0[c} B_{q]a^1} \\
&\qquad\qquad\qquad\qquad + \frac{2}{3} \bg_{ca^0} C_{a^1}{}^{(rs)p} C_{qrsp}
\bigl) \\
&\qquad +X^PX^Q \bigl( -\na_{a^0} B_{a^1c} - 2(n-4) P_{a^0}{}^q A_{q(a^1c)}
+ \frac{4}{3} \bg_{ca^0} C_{a^1}{}^{(rs)p} A_{p(rs)} \bigl) \Bigr] \\
& + X_{B} X_C \ph_{a^1}{}^{PQ}
\end{align*}
for some $\ph_{a^1}{}^{PQ} \in \cE_{a^1}{}^{PQ}$ after some computation. 
Using the last display, it is not difficult to
verify
\begin{align*}
(\pa^*_{BC} d^\na \ol{\Ps})_{a^1BC}{}^{PQ} = 
- \frac{n}{2} \ol{\Ps}_{a^1BC}{}^{PQ} 
- \frac{2}{9}(n-2) (\pa^*_{BC} (\Ps \wedge \Ps))_{a^1BC}{}^{PQ}.
\end{align*}
Thirdly, one easily derives
$\ol{\Ps}_{a^0BC}{}^{RS} {\Ps}_{a^1RS}{}^{PQ} = 
-\frac{n-4}{3}{\Ps}_{a^0BC}{}^{RS} {\Ps}_{a^1RS}{}^{PQ}$.
Hence we finally obtain
\begin{equation} \label{Ps=}
\ol{\ol{\Ps}}_{a^1BC}{}^{PQ} = 
\frac{1}{2} (\pa^*_{BC} \ol{\ol{\Om}})_{a^1BC}{}^{PQ} = 
-\frac{8}{9n}(n-3) (\pa^*_{BC} (\Ps \wedge \Ps))_{a^1BC}{}^{PQ},
\end{equation}
where $-\frac{8}{9n}(n-3) = -\frac{4}{9n}(n-2) - \frac{4}{9n}(n-4)$.

In the last step we need the action of $\Box_1$ on the subbundle
$\cE_{a^1}[-2] \subseteq \cE_{a^1(BC)_0}$ corresponding to the injector
$X_BX_C$. This is the scalar $-(n-1)$, so by adding
$\frac{1}{n-1} \ol{\ol{\Ps}}_{a^1BC}{}^{PQ}$ to $\ol{\ol{\na}}_a$ we obtain
the resulting prolongation covariant derivative
$$
\wt{\na}_a := \na_a +\frac{2}{3} {\Ps}_{aBC}{}^{PQ} 
+ \frac{2}{n} {\ol{\Ps}}_{aBC}{}^{PQ}
+ \frac{1}{n-1} \ol{\ol{\Ps}}_{aBC}{}^{PQ}: \cE_{(PQ)} \to \cE_{a(BC)}.
$$

\begin{proposition}
The prolongation connection
$\widetilde{\na}: \cE_{(BC)} \to \cE_{a(BC)}$
in the conformal geometry has the form
$\widetilde{\na}_a F_{BC}  = \na_a F_{BC} 
+ \frac{2}{3}\wt{\Ps}_{aBC}{}^{PQ} F_{PQ}$
where
\begin{align*}
&\wt{\Ps}_{aBC}{}^{PQ} = 
- Z_{(B}^{\ b}Z_{C)}^{\,c} \Bigl[ 
X^{(P}Z^{Q)q} C_{a(bc)q} + X^PX^Q A_{a(bc)} \Bigr] \\
&\quad + X_{(B}^{}Z_{C)}^{\,c} \Bigl[
-\frac{4}{n} Z^{p(P}Z^{Q)q} C_{a(pq)c} 
+2 X^{(P}Z^{Q)q} \bigl( A_{a(cq)} + \frac{n-4}{n} A_{q(ac)} \bigr)  \\
&\quad\qquad\qquad\quad\
-\frac{1}{n} X^PX^Q B_{ac} \Bigr] \\
&\quad + X_BX_C \Bigl[
-\frac{4}{n} Z^{p(P}Z^{Q)q} A_{a(pq)}
+\frac{1}{n} X^{(P}Z^{Q)q} \bigl( B_{aq} + 
\frac{4}{3(n-1)} C_a{}^{(rs)p} C_{qrsp} \bigr) \\
&\quad\qquad\qquad\quad\
+ \frac{4}{3n(n-1)} X^PX^Q C_a{}^{(rs)p} A_{prs} \Bigr].
\end{align*}
\qed
\end{proposition}
\end{example}

\begin{example}
The prolongation covariant derivative $\wt{\na}$ on tractor
form bundles $\cE_{A^0\form{A}}$, $\form{A} = \form{A}^k$ was computed in
\cite{HaCKF}. Consider a section 
$F_{A^0\form{A}} = \Y_{A^0\form{A}}^{\quad \form{a}} \si_\form{a}
+ \frac{1}{k+1} \Z_{A^0\form{A}}^{\;a^0\form{a}} \mu_{a^0\form{a}} 
+ \W_{A^0\form{A}}^{\quad \dform{a}} \nu_{\dform{a}} 
+ \X_{A^0\form{A}}^{\quad \form{a}} \rh_\form{a} \in \cE_{A^0\form{A}}$.
Then 
\begin{align*}
\wt{\na}_c & F_{A^0\form{A}} = {\na}_c F_{A^0\form{A}}
+ \frac{1}{2} \Z_{A^0\form{A}}^{\;a^0\form{a}} \Bigl[
C_c{}^p{}_{a^0a^1} \si_{p\dform{a}} 
+ \frac{k-1}{n} \bg_{ca^0} C_{a^1a^2}{}^{pq} \si_{pq\ddform{a}} \Bigr] \\
& -\frac{k(k-1)}{2n(n-k)} \W_{A^0\form{A}}^{\quad \dform{a}} \Bigl[
(n-2) C_{ca^2}{}^{pq} \si_{pq\ddform{a}} 
- (k-2) C_{a^2a^3}{}^{pq} \si_{cpq\dddform{a}} \Bigr] \\
&+ \X_{A^0\form{A}}^{\quad \form{a}} \Bigl[ 
- A_c{}^p{}_{a^1} \si_{p\dform{a}}
-\frac{(k-1)(k-2)}{2nk} \bg_{ca^1} C_{a^2a^3}{}^{pq} \nu_{pq\dddform{a}} \\
& \quad + \frac{k-1}{2(n-k)} \Bigl(
\frac{n-2k}{2n} (\na_c C_{a^1a^2}{}^{pq}) \si_{pq\ddform{a}}
+ \bg_{ca^1} A^{pq}{}_{a^2} \si_{pq\ddform{a}} \\
& \qquad\qquad\qquad\quad -2 A_{ca^1}{}^p \si_{p\dform{a}} 
- A_{a^1a^2}{}^p \si_{cp\ddform{a}} + C_{ca^1}{}^{pq} \mu_{pq\dform{a}} \\
& \qquad\qquad\qquad\quad  
+ \frac{n(n-k+1)-2k}{nk} C_c{}^p{}_{a^1a^2} \nu_{p\ddform{a}}
- \frac{k}{n} C_{a^1a^2}{}^{pq} \mu_{cpq\ddform{a}} \Bigl) \Bigr], 
\end{align*}
cf.\ \cite[Remark 4.2]{HaCKF}. 

\vspace{1ex}

The prolongation covariant derivative $\wt{\na}$ simplifies 
for $k=2$ in dimension $n=4$. Then we have (at least locally) 
the conformal volume form 
\begin{equation} \label{vol}
\epsilon_\form{c} \in \cE_\form{c}[4] \quad \text{such that} \quad
\epsilon^\form{c} \epsilon_\form{c} = 4!,  \ \text{i.e.} \ 
\epsilon^\form{e} \epsilon_\form{c} = 4! 
\de_{c^1}^{e^1}\, \de_{c^2}^{e^2}\, \de_{c^3}^{e^3}\, \de_{c^4}^{e^4}, 
\end{equation}
where $\form{c} = \form{c}^4$, $\form{e} = \form{e}^4$. Recall
$\na \epsilon =0$ for any connection $\na$ from the conformal class.
Then the Hodge--star operator 
$\ast: \cE_{\form{a}^k} \to \cE_{\form{a}^{4-k}}$, $k=0,\ldots,4$ has the form 
$(\ast f)_{\form{a}^k} = \epsilon_{\form{a}^k}{}^{\form{r}^{4-k}} 
f_{\form{r}^{4-k}}$. 
The eigenvalues of $\ast$ for $k = 2$ are $\pm 2$. 
The induced tractor volume form 
$E_{\form{C}^6}  = -30 \W_{\form{C}^6}^{\;\form{c}^4} \epsilon_{\form{c}^4} 
\in \cE_{\form{C}^6}$ yields analogously the tractor Hodge--star operator 
$\ast: \cE_{\form{B}^\ell} \to \cE_{\form{B}^{6-\ell}}$. 
The eigenvalues of $E$ for $\ell=3$ are $\pm 6$.

Henceforth we assume $k=2$ and $n=4$ and $*F = 6F$. If not stated otherwise,
all form indices will have valence 2, e.g. $\form{A} = \form{A}^2$ or
$\form{a} = \form{a}^2$. Our normalization of volume forms $E$ and $\epsilon$
means that 
\begin{equation} \label{star}
*\si = 2\si, \quad \ast \mu = -3 \nu, \quad \ast \nu = 2 \mu,
\quad *\rh = -2\rh,
\end{equation}
i.e.\ $\si_\form{a}$ is self-adjoint.
Using this and \nn{vol} one easily verifies
\begin{equation} \label{proj}
\bg_{ca^0} C_{\form{a}}^{}{}^\form{r} \si_\form{r} = 
-2 C_c{}^p{}_\form{a} \si_{pa^0}, \quad
C_\form{a}{}^\form{r} \mu_{c\form{r}} = 
-2 C_{ca^1}{}^\form{r} \mu_{a^2\form{r}}^{}.
\end{equation}
Thus the prolongation covariant derivative $\wt{\na}$ has the form
\begin{align*}
\wt{\na}_c F_{A^0\form{A}} =& \na_c F_{A^0\form{A}}
+ \frac{1}{4} \Z_{A^0\form{A}}^{\;a^0\form{a}} 
C_c{}^p{}_\form{a} \si_{pa^0} 
- \frac{1}{4} \W_{A^0\form{A}}^{\quad a} C_{ca}{}^\form{r} \si_\form{r} \\
&+ \frac{1}{4} \X_{A^0\form{A}}^{\quad \form{a}} \bigl[ 
- 4 A_c{}^p{}_{a^1} \si_{pa^2} + \bg_{ca^1} A^\form{r}{}_{a^2} \si_\form{r} 
- 2 A_{ca^1}{}^p \si_{pa^2} \\
&\qquad\qquad\ - A_\form{a}{}^p \si_{cp} 
+ 2 C_{ca^1}{}^\form{r} \mu_{a^2\form{r}} + C_c{}^p{}_\form{a} \nu_p \bigr], 
\end{align*}

The connection $\wt{\na}$ simplifies considerably for half-flat structures, i.e.\
when
\begin{equation} \label{half}
\epsilon_\form{a}{}^\form{r} C_\form{rb} + 
\epsilon_\form{b}{}^\form{r} C_\form{ar} = 4\la C_\form{ab}, \quad
\la \in \{+1,-1\}.
\end{equation}
The self-adjoint structure $\la=1$ equivalently means 
$C_\form{a}{}^{\form{r}} f_\form{r}=0$ for every anti-self-adjoint two
form $f_\form{a}$ and the anti-self-adjoint structure $\la=-1$ analogously 
means $C_\form{a}{}^{\form{r}} f_\form{r}=0$ for every self-adjoint 
$f_\form{a}$. It follows from \nn{half}, \nn{star} and \nn{vol} that
\begin{equation} \label{equal}
C_c{}^p{}_\form{a} \nu_p = \la C_\form{a}{}^\form{r} \mu_{c\form{r}}^{}.
\end{equation}


We shall discuss the anti-self dual case $\la=-1$ in detail. A short 
computation reveals 
$$
C_\form{a}{}^{\form{r}} \si_\form{r}=0, \quad
A^\form{r}{}_a \si_\form{r} =0 \quad \text{and} \quad
A_\form{a}{}^p \si_{cp} = 2 A_{a^1c}{}^p \si_{a^2p},
$$
where the second and the third equally follow by applying
$\na^{a^1}$ and $\na_{a^0}$, respectively, to the first one and using
$\na_{a^0} C_\form{ar} = 2 \bg_{a^0r^1} A_{\form{a}r^2}$.
(Note the last equality says $A_{[\form{a}}{}^p \si_{c]p}=0$.) 
From the last display and \nn{equal} for $\la=-1$ we finally obtain
the following:

\begin{proposition}
Consider an anti-self-dual conformal structure in the dimension 4.
Then the prolongation connection
$\widetilde{\na}: \cE^+_{[A^0\form{A}]} \to \cE^+_{c[A^0\form{A}]}$,
$\form{A} = \form{A}^2$ on
the bundle of self-dual tractor 3-forms 
$\cE^+_{[A^0\form{A}]} \subseteq \cE_{[A^0\form{A}]}$ has the form
\begin{equation*}
\wt{\na}_c F_{A^0\form{A}} = \na_c F_{A^0\form{A}}
+ \X_{A^0\form{A}}^{\quad \form{a}} \bigl[ 
- 2 A_{c(pa^1)} \si^p{}_{a^2} 
+ \frac{1}{2} C_c{}^p{}_\form{a} \nu_p \bigr].
\end{equation*}
for $F_{A^0\form{A}} \in \cE^+_{[A^0\form{A}]}$ where
$\si_\form{a} = 3\X^{A^0\form{A}}_{\quad \form{a}} F_{A^0\form{A}}$ and 
$\nu_{a} = -6\W^{A^0\form{A}}_{\quad a} F_{A^0\form{A}}$. 
\end{proposition}

Note a modification of $\na$ on $\cE^+_{A^0\form{A}}$ was also obtained
in \cite[(2.27)]{DuTo4dim} where the spinorial notation is used. 
%
\end{example}

\comments{
\section{Commutation of higher projective BGG squares for the tractor forms} 

Let us denote by $\na$ the normal tractor connection. We shall use the notation 
from Section \ref{forms} for forms. In particular, all sequentially labeled 
indexes are implicitly skewed over.

We denote by $\cT$ the standard projective tractor bundle. The composition 
series of $(\cT^*)^k := \bigwedge^k\cT^*$ is 
$\cE_{\dform{a}}(k) \lpl \cE_{\form{a}}(?)$
and the normal tractor covariant derivative $\na$ and its curvature 
$\Om$ take the form
$$
\na_c
\begin{pmatrix}
\si_{\dform{a}} \\ \mu_{\form{a}}
\end{pmatrix} = 
\begin{pmatrix}
\na_c \si_{\dform{a}} + \mu_{c\dform{a}} \\ 
\na_c \mu_{\form{a}} - k P_{ca^1} \si_{\dform{a}}
\end{pmatrix}, 
\quad
\Om_{c^0c^1} \sharp \begin{pmatrix}
\si_{\dform{a}} \\ \mu_{\form{a}}
\end{pmatrix} = 
\begin{pmatrix}
(k-1)W_{c^0c^1}{}^p{}_{a^2} \si_{p\ddform{a}} \\ 
kW_{c^0c^1}{}^p{}_{a^1}\mu_{p\dform{a}} + \ast
\end{pmatrix}.
$$
Here $\form{a} = \form{a}^k$.
We shall consider study sections
$f_{\form{c}} \in \cE_{\form{c}} \otimes \Ga((\cT^*)^k)$, i.e.\ tractor 
indices are supressed in the notation for $f$. Here $\form{c} = \form{c}^\ell$.

Following Theorem \ref{?}, we shall start with the sequence $E_i = d^\na$. 
Our aim is to find a suitable modification 
$D: \cE_\form{c} \otimes (\cT^*)^k \to 
\cE_{c^0\form{c}} \otimes (\cT^*)^k$ of $d^\na$ such that
$(D \circ D f)_{c^{-1}c^0\form{c}}) \in 
\Ker \partial^* \subseteq \cE_{c^{-1}c^0\form{c}}$ for every $f_\form{c} \in 
\cE_\form{c} \otimes (\cT^*)^k$. Specifically, we shall use
certain (algebraic) operator 
${\Ps}_{c^0} \in \cE_c \otimes \End((\cT^*)^k)$
and put $(\Ps(f))_{c^0\form{c}} := ({\Ps} \wedge f)_{c^0\form{c}} \in 
\Im \partial^* \subseteq \cE_{c^0\form{c}} \otimes \End((\cT^*)^k)$,
and we put $(Df)_{c^0\form{c}} := (d^\na f)_{c^0\form{c}} +
(\Box_k)^{-1} {\Ps}_{c^0}f_\form{c} \in 
\cE_{c^0} \otimes (\cT^*)^k$. We shall usually write the endomorphism 
$\Ps$ as $\Ps_{c^0} f_\form{c}$ as the notation already requires the skew 
symmetrization. Here $\Box_k$ is a scalar multiple determined by the Kostant's 
Laplacian.

Now we shall describe the difference between $d^\na$ and $D$ on 
$\cE_\form{c} \otimes \End((\cT^*)^k)$. Consider the section
$$
f_\form{c} =  
\begin{pmatrix}
\si_{\dform{a}} \\ \mu_{\form{a}}
\end{pmatrix}
\in \cE_{\form{c}} \otimes \End((\cT^*)^k)
$$
and 
$$
(d^\na d^\na f)_{c^{-1}c^0\form{c}} = \Om_{c^{-1}c^0} \sharp f_\form{c} = 
\begin{pmatrix}
W_{c^{-1}c^0}{}^p{}_{a^2} \si_{\form{c}p\ddform{a}} \\ \ast
\end{pmatrix}
\in \cE_{c^{-1}c^0\form{c}} \otimes \End((\cT^*)^k).
$$
Now we need $\partial^*$ of the previous display and this will define the 
endomorphism $\Ps$. We need it only up to a scalar multiple at the moment
and a short computation shows
$$
(\Ps(f))_\form{c^0\form{c}} = (\partial^* \Om \sharp f)_{c^0\form{c}} = 
\begin{pmatrix}
0 \\
W_{a^1a^2}{}^p{}_{c^0} \si_{\form{c}p\ddform{a}} 
- \ell W_{c^{-1}c^0}{}^p{}_{a^1} \si_{a^2\dform{c}p\ddform{a}} 
\end{pmatrix}.
$$

We use the modification $D := d^\na + \al \Box_k^{-1}\Ps$ and it remains 
to determine the operator $\Box_k$, i.e.\ how we need to rescale particular 
($\frak{g}_0$-)irreducible components of $\Ps(f)$.

\section{More complicated examples: projective geometry}

The standard projective tractor bundle $\cE^A$ has the composition
series $\cE^A = \cE^a(1) \lpl \cE(?)$. We shall write section of this bundle
as 
$$
f^A = 
\begin{pmatrix} 
\si^a \\ \rh
\end{pmatrix}
= Y^A_a \si^a + X^A \rh
$$
where $X_A: \cE(?) \to \cE^A$ is invariant. 
The covariant derivative is then
$$
\na_c F^A = 
\begin{pmatrix} 
\na_c \si^a \\ \na_c \rh
\end{pmatrix}, 
\quad \text{i.e.} \quad
\na_c Y^A = , \na_c X^A =
$$
}

\section{Almost Grassmannian geometry}

A complex almost Grassmannian (or AG--) structure on a smooth manifold 
$M$ is given 
by two auxiliary vector bundles $\cE^A$ and $\cE_{A'}$ and the identification
\begin{equation} \label{iso}
\cE^a = \cE_{A'} \otimes \cE^A = \cE^A_{A'}, \quad
\bigwedge\!\!{}^q \cE^A \cong \bigwedge\!\!{}^p \cE_{A'},
\end{equation}
where $p$ is the rank of $\cE_{A'}$ and $q$ is the rank of $\cE^A$.
In fact, all results we obtain hold for all real forms of a given complex 
geometry, \cite{GS}. Motivated by the case $p=q=2$ when the structure is the spin
conformal structure, we shall term $\cE^A$ and $\cE_{A'}$ spinor bundles. 

Following \cite{GS} and equation \nn{iso}, we adopt the convention
$$
\cE[-1] \cong \cE_{\form{A}^q} \cong \cE^{\form{B'}^p}, \quad
\cE[1] \cong \cE^{\form{A}^q} \cong \cE_{\form{B'}^p}
$$
for line bundles. This isomorphism is given explicitly
by the tautological section $\vol_{\form{A}^q} \in \cE_{\form{A}^q}[1]$
as $\cE[-1] \ni f \mapsto f\vol_{\form{A}^q} \in \cE_{\form{A}^q}$.
A choice of a \idx{scale} $\xi \in \cE[1]$ is equivalent to the choice of spinor
volume forms $\vol^\xi_{\form{A}^q} := 
\xi^{-1} \vol_{\form{A}^q} \in \cE_{\form{A}^q}$, and analogously for 
$\cE^{\form{A'}^p}$. 

Our convention for the torsion $T_{ab}{}^c$ and the curvature $R_{ab}{}^d{}_c$ 
of a covariant derivative $\na_a$ on $TM$ are given by the equation
$$
2\na_{[a} \na_{b]}v^c = T_{ab}{}^d \na_d v^c + R_{ab}{}^c{}_d v^d.
$$
Summarizing \cite[Theorem 2.1]{GS}, for a scale $\xi \in \cE[1]$ on an
AG--structure there are unique covariant derivatives on $\cE^A$ and $\cE_{A'}$
such that the torsion $F_A^{A'}{}_B^{B'}{}_{C'}^C$ of the induced covariant derivative on 
$TM$ is totally trace-free, the induced covariant derivative preserves 
\nn{iso} and in addition, $\xi$ is parallel. We denote this class
of covariant derivatives, parametrized by sections of $\cE[1]$, by $[\na]$.
Changing the scale $\xi \to \hat{\xi} = e^{\Up}\xi \in \cE[1]$ with $\Up$
a smooth function, the covariant derivative $\na$ changes to $\hat{\na}$
in a way that
\begin{align} \label{Gtrans}
\begin{split} 
&\hat{\na}_A^{A'} u^C \; = \na_A^{A'} u^C + \de_A^C \Up_B^{A'} u^B,
\quad \mbox{for}\ u^C \in \cE^A, \\
&\hat{\na}_A^{A'} u_{C'} = \na_A^{A'} u_{C'} + \de_{C'}^{A'} \Up_A^{B'} u_{B'},
\quad \mbox{for}\ u_{C'} \in \cE_{C'}, \\
&\hat{\na}_A^{A'} v_B \; = \na_A^{A'} v_B - \Up_B^{A'} v_A,
\quad \mbox{for}\ v_B \in \cE_B. \\
&\hat{\na}_A^{A'} v^{B'} = \na_A^{A'} v^{B'} - \Up_A^{B'} v^{A'},
\quad \mbox{for}\ v_B \in \cE_B \ \ \mbox{and also} \\
&\hat{\na}_a f = \na_a f +w \Up_a f,
\quad \mbox{for}\ f \in \cE[w] 
\end{split}
\end{align}
where $\Up_a = \na_a \Up$. We use hat sign 
to denote quantities corresponding to the changed scale 
$\hat{\xi} = e^{\Up}\xi$ from now on without further notice. 

Given $\na \in [\na]$, we denote all covariant derivatives on tensor products 
of $\cE^A$ ans $\cE_{A'}$ also by $\na$. The curvature on spinor bundles is 
given by
$$
(2\na_{[a}\na_{b]} - T_{ab}{}^d \na_d) v^C = R_{ab}^{}{}_D^C v^D,
\quad
(2\na_{[a}\na_{b]} - T_{ab}{}^d \na_d) v_{D'} = -R_{ab}^{}{}_{D'}^{C'} v_{C'}.
$$
The curvature of $\na$ is $R_{ab}^{}{}_c^d = 
R_{ab}^{}{}_{C}^{D} \de_{D'}^{C'} - R_{ab}^{}{}_{D'}^{C'} \de_{C}^{D}$,
where $R_{ab}^{}{}_{C'}^{D'}$ and $R_{ab}^{}{}_D^C$ are trace--free on the 
spinor indices displayed. The relations
\begin{align*}
&R_{ab}^{}{}_D^C = U_{ab}^{}{}_D^C 
- \de_B^C \Rho_A^{A'}{}_D^{B'} + \de_A^C \Rho_B^{B'}{}_D^{A'}, \\
&R_{ab}^{}{}_{D'}^{C'} = U_{ab}^{}{}_{D'}^{C'} 
+ \de_{D'}^{B'} \Rho_A^{A'}{}_B^{C'} - \de_{D'}^{A'} \Rho_B^{B'}{}_A^{C'}, 
\end{align*}
together with the condition
$U_R^{A'}{}_B^{B'}{}_A^R - U_A^{R'}{}_B^{B'}{}_{R'}^{A'}=0$
(and the algebraic Bianchi identity)
determine $U_{ab}^{}{}_D^C$, $U_{ab}^{}{}_{D'}^{C'}$ and 
the \idx{Rho--tensor} $\Rho_{ab}$. In more details, the curvature
on the (co)tangent bundle is
$$
R_{ab}^{}{}_d^c = U_{ab}^{}{}_d^c +
\de_{C'}^{D'} \de_A^C \Rho_B^{B'}{}_D^{A'} 
- \de_{C'}^{D'} \de_B^C \Rho_A^{A'}{}_D^{B'}
+ \de_D^C \de_{C'}^{A'} \Rho_B^{B'}{}_A^{D'}
- \de_D^C \de_{C'}^{B'} \Rho_A^{A'}{}_B^{D'}
$$
where $U_{ab}^{}{}_d^c = U_{ab}^{}{}_D^C \de_{C'}^{D'} 
- U_{ab}^{}{}_{C'}^{D'} \de_D^C$. In this form, tensors $U$ are determined by 
$U_{rb}^{}{}_a^r = U_R^{A'}{}_B^{B'}{}_A^R - U_A^{R'}{}_B^{B'}{}_{R'}^{A'}=0$.
(Note the previous display means the decomposition $U = R + \pa \Rho$ where 
$U$ is $\pa^*$-closed, cf.\ the theory of Weyl structures in \cite{CSbook}.)
Furthermore,
\begin{equation} \label{RhoT}
U_{ab}^{}{}_C^C = -U_{ab}^{}{}_{C'}^{C'} = 2\Rho_{[ab]} 
\quad \mbox{and} \quad
-2(p+q)\Rho_{[ab]} = \na_c T_{ab}{}^c
\end{equation}
where the last identity follows from the algebraic Bianchi identity.

We will be mostly interested in the case $p=2$ and $q>2$. In this case,
the only invariants are the trace--free part of 
$T_{\;[A\,B]}^{(A'\!B')}{}_{C'}^C$ and the trace--free part 
of $U_{\; (A\,B\,C)}^{[A'B']D}$, \cite{GS}. That is, if these two vanish, 
the geometry is locally isomorphic to the homogenous model.
Finally note that using the algebraic Bianchi identity 
we obtain
\begin{align} \label{p=2Tr}
\begin{split}
&U_{\; (A\,B)R'}^{R'[A'B']} = U_{\ R\,(A\;B)}^{[A'B']R}
= U_{\; [A\,B]R'}^{R'(A'B')} = U_{\ R\,[A\;B]}^{(A'B')R'}=0, \\ 
&U_{\; (A\,B)R'}^{R'(A'B')} = U_{\ R\,(A\;B)}^{(A'B')R}
= \frac{1}{q} T_r^{}{}_{(A}^{(A'|e|} T_{\:B\,)\,e}^{B')\ r}  \\
&U_{\; [A\,B]R'}^{R'[A'B']} = U_{\ R\,[A\;B]}^{[A'B']R}
= -\frac{1}{q+4} T_r^{}{}_{[A}^{[A'|e|} T_{\:B\,]\,e}^{B']\ r}.
\end{split}
\end{align}

\subsection{Grassmannian tractor calculus.} 
We follow \cite{GS} here. The standard tractor bundle
is the (spinor tractor) bundle $\cE^\al = \cE^A \lpl \cE^{A'}$ and we denote
its dual by $\cE_\al  = \cE_{A'} \lpl \cE_A$. (That is, we use Greek letters 
for spinor tractor abstract indices.) 
Using the injectors $Y^\al_A \in \cE^\al_A$, $X^\al_{A'} \in \cE^\al_{A'}$
and $Y_\al^{A'} \in \cE_\al^{A'}$, $X_\al^{A} \in \cE_\al^{A}$,
sections of $\cE^\al$ and $\cE_\al$ are written conveniently as
$$
\begin{pmatrix} 
\si^A \\ \rh^{A'}
\end{pmatrix}
= Y^\al_A \si^A + X^\al_{A'} \rh^{A'} \in \cE^\al,
\quad \text{resp.} \quad
\begin{pmatrix} 
\nu_{A'} \\ \mu_A
\end{pmatrix}
= Y_\al^{A'} \nu_{A'} + X_\al^{A} \mu_A \in \cE_\al.
$$
Splittings of $\cE^\al$ and $\cE_\al$ are parametrised by choice of
the scale $\xi \in \cE[1]$. The change of the splitting has the form
\begin{align*}
&\widehat{
\begin{pmatrix} 
\si^A \\ \rh^{A'}
\end{pmatrix}
} = 
\begin{pmatrix} 
\si^A \\ \rh^{A'} - \Up_B^{A'} \si^B  
\end{pmatrix},
\ \ \mbox{i.e.} \ \
\hat{Y}^\al_A = Y^\al_A + X^\al_{B'} \Up_A^{B'},\ 
\hat{X}^\al_{A'} = X^\al_{A'}
\quad \text{and} \\
&\widehat{
\begin{pmatrix} 
\nu_{A'} \\ \mu_A
\end{pmatrix}
} = 
\begin{pmatrix} 
\nu_{A'} \\ \mu_A + \Up_A^{A'} \nu_{A'}
\end{pmatrix},
\ \ \mbox{i.e.} \ \
\hat{Y}_\al^{A'} = Y_\al^{A'} - X_\al^{B} \Up_B^{A'},\ 
\hat{X}_\al^{A} = X_\al^{A}
\end{align*}
That is, the sections $X^\al_{A'}$ and $X_\al^{A}$ are invariant
and $Y^\al_A$ and $Y_\al^{A'}$ depend on the choice
of the scale. 
They are normalized in such a way that 
$Y^\be_B X_\al^B + Y_\al^{B'} X^\be_{B'} = \de_\al{}^\be$, i.e.\ 
$X_\al^B Y^\al_A =\de_{A}{}^{B}$ and 
$X^\al_{A'} Y_\al^{B'} = \de_{A'}{}^{B'}$. 

The normal covariant tractor derivative is given by
\begin{align*}
&\na_A^{P'} \!\!
\begin{pmatrix} 
\si^B \\ \rh^{B'}
\end{pmatrix} = 
\begin{pmatrix} 
\na_A^{P'} \si^B + \rh^{P'} \de_A{}^B \\ 
\na_A^{P'} \rh^{B'} - \Rho_A^{P'\!}{}_B^{B'} \si^B
\end{pmatrix} 
\ \text{and} \
\na_A^{P'} \!\!
\begin{pmatrix} 
\nu_{B'} \\ \mu_B
\end{pmatrix} = 
\begin{pmatrix} 
\na_A^{P'} \nu_{B'} - \de_{B'}^{P'} \mu_A  \\ 
\na_A^{P'} \mu_B + \Rho_A^{P'\!}{}_B^{B'} \nu_{B'}
\end{pmatrix}.
\end{align*}
That is, 
\begin{align*}
&\na_A^{P'} Y^\al_B = - X^\al_{B'} \Rho_A^{P'\!}{}_B^{B'}, 
\na_A^{P'} X^\al_{B'} = Y^\al_A \de_{B'}^{P'}
\ \text{and} \ \\
&\na_A^{P'} Y_\al^{B'} = X_\al^{B} \Rho_A^{P'\!}{}_B^{B'},\ 
\na_A^{P'} X_\al^{B} = -Y_\al^{P'} \de_A^B.
\end{align*}
Its curvature $\Om_{ab}^{}{}^\al_\be$ is trace-free on the spinor tractor 
bundle and has the explicit form
\begin{align*}
\Om_{ab}^{}{}^\al_\be = & - Y^\al_C Y_\be^{C'} T_{ab}^{}{}^C_{C'} 
+ Y^\al_C X_\be^{D} U_{ab}^{}{}^C_{D} 
+ X^\al_{C'} Y_\be^{D'} U_{ab}^{}{}^{C'}_{D'} \\
&+ X^\al_{C'} X_\be^{C} Q_{ab}^{}{}^{C'}_C \in
\cE_{[ab]}^{}{}^\al_\be \subseteq \cE_{[ab]} \otimes 
\text{trace-free} (\cE^\al_\be)
\end{align*}
where 
$Q_{abc} = -2 \na_{[a}\Rho_{b]c} + T_{ab}{}^e \Rho_{ec} \in \cE_{[ab]c}$
and $\text{trace-free} (\cE^\al_\be) = \cA$ is the adjoint tractor bundle.
That is, $(\na_a\na_b - \na_b\na_a - T_{ab}{}^e\na_e) f^\al  = 
\Om_{ab}^{}{}^\al_\be f^\be = (\Om \sharp f)_{ab}{}^\al = 
\Om_{ab} \sharp f^\al$ in our notation. 

\vspace{1ex}
The inclusions $\io$ and $\bar{\io}$ from \ref{KoLap} are of the form
$Y^\al_{A^0} Y_\be^{A^0{}'}: \cE_\form{a} \stackrel{\bar{\io}}{\to} 
\cE_{a^0\form{a}}{}^\al{}_\be$ and 
$X^\al_{{A^1}'} X_\be^{{A^1}}: \cE_\form{a} \stackrel{\io}{\to} 
\cE_{\dform{a}}{}^\al{}_\be$, where we use the identification
$\cE_{a^0}^{} = \cE_{A^0}^{A^0{}'}$ and $\cE^{a^1} = \cE_{A^1{}'}^{A^1}$.
Therefore
\begin{align*}
&\pa: \cE_\form{a} \otimes \cT \ni f_\form{a} \mapsto
Y^\al_{A^0} Y_\be^{A^0{}'} f_\form{a} \stackrel{\sharp}{\longrightarrow} 
\cE_{a^0\form{a}} \otimes \cT \quad \text{and} \\
&\pa^*: \cE_\form{a} \otimes \cT \ni f_\form{a} \mapsto
X^\al_{{A^1}'} X_\be^{{A^1}} f_\form{a} \stackrel{\sharp}{\longrightarrow} 
\cE_{\dform{a}} \otimes \cT
\end{align*}
for any subbundle $\cT$ of 
$\bigotimes \cE_\al \otimes \bigotimes \cE^{\be} \otimes \cE[w]$. This does
not cover all tractor bundles but will be sufficient in the examples treated 
below.

\vspace{1ex}

Henceforth we assume $p=2$, $q>2$. Note we have the decomposition
$\Om_{ab}^{}{}^\al_\be = \Om_{\,(AB)}^{[A'B']}{}^\al_\be + 
\Om_{\,[A\,B]}^{(A'B')}{}^\al_\be$, where the component 
$\Om_{\,[A\,B]}^{(A'B')}{}^\al_\be$ vanishes in the torsion--free case.

\subsection{Skew symmetric tractors and tractor forms}
We shall also need tractor bundles 
$\bigwedge^\ell \cE^\al = \cE^{\boldsymbol{\al}}$ with the notation for the multiindex ${\boldsymbol{\al}} = {\boldsymbol{\al}}^\ell$. Since
$\bigwedge^\ell \cE^\al \cong \bigwedge^{q+2-\ell} \cE_\be$ (we assume
orientability here), these are just tractor forms. Specifically, the 
case $\ell = q+1$ is just the bundle $\cE_\be$. 

It follows from the structure of $\cE^\al$ that
$$
\cE^{\boldsymbol{\al}} =
\cE^\form{A} \lpl \cE^{B'\dform{A}} \lpl \cE^{[B'C']\ddform{A}},
\qquad \boldsymbol{\al} = \boldsymbol{\al}^\ell, \ 
\form{A} = \form{A}^\ell, \ 2 \leq \ell \leq q.
$$
Of course we have the isomorphism
$\cE^{[B'C']\ddform{A}} \cong \cE^{\ddform{A}}[-1]$ using the spinor volume
form $\vol_{B'C'} \in \cE_{[B'C']}[-1]$ but it turns out more convenient for 
the computation to use the form as in the display.

We put
\begin{align*}
&\Y^{\boldsymbol{\al}}_{\form{A}} 
= Y^{[\al^1}_{\ A^1} \ldots Y^{\al^\ell]}_{\,A^\ell}
\in \cE^{\boldsymbol{\al}}_{\,\form{A}}, \quad
\W^{\,\boldsymbol{\al}}_{B'\dform{A}} 
= X^{[\al^1}_{\; B'} Y^{\al^2}_{\ A^2} \ldots Y^{\al^\ell]}_{\,A^\ell}
\in \cE^{\,\boldsymbol{\al}}_{\,B'\dform{A}}, \\
&\X^{\,\boldsymbol{\al}}_{B'C'\dform{A}} 
= X^{[\al^1}_{\; B'} X^{\al^2}_{\ C'} Y^{\al^3}_{\ A^3} 
\ldots Y^{\al^\ell]}_{\,A^\ell}
\in \cE^{\,\boldsymbol{\al}}_{[B'C']\ddform{A}}
\end{align*}
where $\X^{\,\boldsymbol{\al}}_{B'C'\dform{A}}$ is invariant and 
$\Y^{\boldsymbol{\al}}_{\form{A}}$ and $\W^{\,\boldsymbol{\al}}_{B'\dform{A}}$
are scale dependent. Finally, the normal tractor connection on these section
is
\begin{align*}
&\na_c \Y^{\boldsymbol{\al}}_{\form{A}} = -\ell\, 
\W^{\,\boldsymbol{\al}}_{B'[\dform{A}} \Rho_{|c|}^{}{}_{A^1]}^{B'}, \quad \\
&\na_c^{} \W^{\,\boldsymbol{\al}}_{B'\dform{A}} = 
\Y^{\boldsymbol{\al}}_{C\dform{A}} \, \de_{B'}^{C'}
- (\ell-1)\, \X^{\,\boldsymbol{\al}}_{B'D'[\dform{A}} 
\Rho_{|c|}^{}{}_{A^2]}^{D'}, \quad \text{and} \\
&\na_c^{} \X^{\,\boldsymbol{\al}}_{B'D'\dform{A}}  = 2\,
\W^{\,\boldsymbol{\al}}_{B'C\ddform{A}} \, \de_{D'}^{C'}. 
\end{align*}


\begin{example}
We shall demonstrate the prolongation covariant derivative $\wt{\na}$ for 
AG--geometries on tractor bundles corresponding to fundamental 
representations. These are bundles $\bigwedge^\ell \cE^{\al}$ for 
$1 \leq \ell \leq q+1$. Since the computation is getting very technical for
$1<\ell<q+1$, we later restrict to torsion--free manifolds.

First we discuss the cases $\cE^\al$ and 
$\cE_\be \cong \bigwedge^{q+1} \cE^{\al}$. Considering $F^\al \in \cE^\al$
and $G_\be \in \cE_\be$, a short computation gives
\begin{align*}
&(\pa^* d^\na d^\na F)_c{}^\al = 
\frac12 X^\al_{D'} X_{\,\om}^{D} S_{\,C\,D}^{D'C'} F^\om 
\quad \text{and} \\
&(\pa^* d^\na d^\na G)_{c\be} = 
- \frac12 X^\om_{D'} X_{\,\be}^{D} S_{\,C\,D}^{D'C'} G_\om^{},
\end{align*}
where
$$
S_{\,C\,D}^{D'C'} = U_{\,A\;B\;R'}^{R'A'B'} = U_{\,R\;A\;B}^{A'B'R}
= \frac{1}{q} T_r^{}{}_{(A}^{(A'|e|} T_{\:B\,)\,e}^{B')\ r} 
-\frac{1}{q+4} T_r^{}{}_{[A}^{[A'|e|} T_{\:B\,]\,e}^{B']\ r}.
$$
Hence we need the action of the Kostant-Laplace operator $\Box$ on
$\cE_{\,C}^{D'C'} = \cE_{\ C}^{(D'C')} \oplus \cE_{\ C}^{[D'C']}$.
The eigenvalues are, respectively, $\frac12(q-1)$ and $\frac12(q+1)$. Therefore
the prolongation connection $\wt{\na}$ has the form
\begin{align*}
&\wt{\na}_c F^\al = \na_c F^\al - X^\al_{\,D'} X_{\,\om}^{D} 
\Bigl[ \frac{1}{q-1} S_{\ C\:D}^{(D'C')} 
+ \frac{1}{q+1} S_{\ C\:D}^{[D'C']} \Bigr] F^\om
\quad \text{for}\  F^\al \in \cE^\al, \\
&\wt{\na}_c G_\be = \na_c G_\be + X^\om_{\,D'} X_{\,\be}^{D} 
\Bigl[ \frac{1}{q-1} S_{\ C\:D}^{(D'C')} 
+ \frac{1}{q+1} S_{\ C\:D}^{[D'C']} \Bigr] G_\om
\quad \text{for}\  G_\be \in \cE_\be.
\end{align*}

It remains to consider the bundles $\cE^{\boldsymbol{\al}}$, 
$\boldsymbol{\al} = \boldsymbol{\al}^\ell$ for $2 \leq \ell \leq q$.
Consider the section
$F^{\boldsymbol{\al}} = \Y^{\boldsymbol{\al}}_{\form{A}} \si^\form{A}
+ \W^{\,\boldsymbol{\al}}_{B'\dform{A}} \mu^{B'\dform{A}}
+ \X^{\,\boldsymbol{\al}}_{B'C'\dform{A}} \rh^{B'C'\ddform{A}}$,
where $\si^\form{A} \in \cE^\form{A}$, 
$\mu^{B'\dform{A}} \in \cE^{B'\dform{A}}$ and 
$\rh^{B'C'\ddform{A}} \in \cE^{[B'C']\ddform{A}}$.
A straightforward computation shows that 
\begin{align*}
(d^\na d^\na  F&)_{de}^{\ \ \boldsymbol{\al}} =
\frac{1}{2} \Om_{de}^{} \sharp F^{\boldsymbol{\al}} = \frac{1}{2} \ell\,
\Om_{de}^{}{}_{\; \om}^{[\al^1} F^{|\om|\boldsymbol{\dform{\al}}]}= \\
= \frac{1}{2} \biggl\{
&\Y^{\boldsymbol{\al}}_{\form{A}} \Bigl[
\ell U_{de}^{}{}_{\; Q}^{[A^1} \si_{}^{|Q|\dform{A}]}
- T_{de}^{}{}_{\; Q'}^{[A^1} \mu_{}^{|Q'|\dform{A}]} \Bigr] \\
&+ \W^{\,\boldsymbol{\al}}_{B'\dform{A}} \Bigl[
(\ell-1) U_{de}^{}{}_{\; Q}^{[A^2} \mu_{}^{|B'Q|\ddform{A}]}
+ \ell Q_{de}^{}{}_{Q}^{B'} \si_{}^{Q\dform{A}}
+ U_{de}^{}{}_{\; Q'}^{B'} \mu_{}^{Q'\dform{A}} \\
&\qquad\qquad 
-2 T_{de}^{}{}_{\; Q'}^{[A^2} \rh_{}^{|B'Q'|\ddform{A}]} \Bigr] 
+ \X^{\,\boldsymbol{\al}}_{B'C'\dform{A}} \ph_{}^{B'C'\ddform{A}}
\biggr\}
\end{align*}
for a section $\ph_{}^{B'C'\ddform{A}} \in \cE_{}^{B'C'\ddform{A}}$.
We need to compute $\pa^*$ of the previous display. 

It turns out the computation
is getting too technical in general, so we compute $\wt{\na}$ in the 
torsion--free case only. That is, we assume $T_{ef}^{}{}_{C}^{C'}=0$
(hence also $S_{\,C\,D}^{D'C'}=0$) from now on. Then we obtain
\begin{align*}
&(\pa^* d^\na d^\na F)_e{}^{\boldsymbol{\al}} =\\
&= \frac{1}{2} (\ell-1) \, \Bigl\{
\ell\, \W^{\,\boldsymbol{\al}}_{B'\dform{A}} 
U_{e}^{}{}_{R\,\ Q}^{B'[A^2} \si_{}^{|QR|\ddform{A}]} \\
&\qquad
+ \X^{\,\boldsymbol{\al}}_{B'C'\ddform{A}} \bigl[
(\ell-2) U_{e}^{}{}_{R\,\ Q}^{C'[A^3} \mu_{}^{|B'QR|\dddform{A}]}
- \ell \,Q_{e}^{}{}_{R\,Q}^{C'B'} \si_{}^{QR\dform{A}} 
- U_{e}^{}{}_{R\,Q'}^{C'B'} \mu_{}^{Q'R\ddform{A}}
\bigr] \Bigr\}.
\end{align*} 
Since $U_{\; A\,B\,C}^{A'B'D} = U_{\; (A\,B\,C)}^{[A'B']D}$ in the 
torsion--free case, we conclude that 
$(\pa^* d^\na d^\na F)_e{}^{\boldsymbol{\al}}=0$. This yields
the surprising result $\wt{\na} = \na$ on $\cE^{\boldsymbol{\al}}$.
The same is obviously true also for $\ell=1$ and $\ell=q+1$. Hence we obtain
\end{example}

\begin{proposition}
The prolongation connection 
$\wt{\na}_c: \cE^{\boldsymbol{\al}} \to \cE_c^{\ \boldsymbol{\al}}$,
$\boldsymbol{\al} = \boldsymbol{\al}^\ell$ for $1 \leq \ell \leq q+1$
on torsion--free AG-manifolds is equal to the normal tractor connection, i.e.\
$\wt{\na} = \na$.
\end{proposition}

\section{Acknowledgement.}
The first author was supported by a Junior Research
Fellowship of The Erwin Schr\"{o}dinger International Institute for
Mathematical Physics
and project P23244-N13 of the "Fonds zur
F\"{o}rderung der wissenschaftlichen Forschung" (FWF).
The second and the third authors were supported by the 
institutional grant MSM 0021620839 and
by the grant GA CR 201/08/397. Most of the work of the 
last author was undertaken while his stay at the 
Max-Planck-Institute f\"{u}r Mathematik in Bonn. The support of MPIM is 
gratefully acknowladged.

\end{document}